\documentclass{article}
\usepackage{amsfonts}
\usepackage{amsmath}
\usepackage{amssymb}
\usepackage{wasysym}
\usepackage{yhmath}
\usepackage{graphicx}

\newtheorem{theorem}{Theorem}
\newtheorem{definition}[theorem]{Definition}
\newtheorem{corollary}[theorem]{Corollary}
\newtheorem{proposition}[theorem]{Proposition}
\newtheorem{example}[theorem]{Example}
\newtheorem{examples}[theorem]{Examples}
\newtheorem{remark}[theorem]{Remark}
\newtheorem{remarks}[theorem]{Remarks}

\newtheorem{lemma}[theorem]{Lemma}
\newtheorem{properties}[theorem]{Properties}

\newtheorem{propmono}[theorem]{Monotonicity of the Busemann cocycle}

\newtheorem{criterium}[theorem]{Criterium}
\newtheorem{approlemma}[theorem]{Uniform Approximation Lemma}

\newtheorem{excesslemma}[theorem]{Excess Lemma}

\def \N{\mathbb{N}}
\def \Q{\mathbb{Q}}
\def \H{\mathbb{H}}
\def \Z{\mathbb{Z}}
\def \R{\mathbb{R}}
\def \Hy{\mathbf{H}}

\def \G{\Gamma}
  
\def \qed{$\Box$}

\begin{document}

\title{\Large On the horoboundary  and the geometry of rays \\ of negatively  curved manifolds.}

\date{\normalsize July, 2010}

\author{\normalsize{Fran\c{c}oise Dal'bo\footnote{Fran\c{c}oise Dal'bo, IRMAR, UniversitŽ de Rennes-I\@, Campus de Beaulieu, 35042 Rennes Cedex -- francoise.dalbo@univ-rennes1.fr},
Marc Peign\'e\footnote{Marc Peign\'e, UMR 6083, Facult\'e des Sciences et Techniques, Parc de Grandmont, 37200 Tours -- peigne@univ-tours.fr},
Andrea Sambusetti\footnote{Andrea Sambusetti, Istituto di Matematica G. Castelnuovo Universit\`a ``La Sapienza'' di Roma      P.le Aldo Moro 5 - 00185 Roma -- sambuset@mat.uniroma1.it}}}

\maketitle

\section{Introduction}

${}$
\vspace{-6mm}

\noindent The problem of understanding the geometry and dynamics of geodesics and rays (i.e. distance-minimizing half-geodesics) on Riemannian manifolds dates back at least to Hadamard \cite{hadamard}, who started to study the qualitative behaviour of geodesics on nonpositevely curved surfaces of $\mathbb{R}^3$. In particular, he first distinguished between different kinds of ends on such surfaces, and introduced the notion of asymptote, which we shall be concerned about in this paper.
\\
Half a century later, in his seminal book \cite{busemann}, Busemann introduced an amazingly simple notion for measuring the   ``angle at infinity''  between rays (now known as the  {\em Busemann function}) as a tool to develop a theory of parallels on geodesic spaces. The Busemann function of a ray $\alpha$   is the two-variables function
\vspace{-6mm}

$$ B_{\alpha}(x,y)= \lim_{t \rightarrow +\infty} d(x, \alpha(t)) - d(\alpha(t), y)$$

\vspace{-1mm}

\noindent and played an important  role (far beyond the purposes of his creator) in the study of complete noncompact Riemannian manifolds. \\
 It  has been used to derive fundamental results in nonnegative curvature such as Cheeger-Gromoll-Meyer's Soul Theorem or  Toponogov' Splitting Theorem (\cite{shiohama}), in the function  theory of harmonic and noncompact symmetric spaces (\cite{anderson}, \cite{macpherson}), and has a special place in the geometry of Hadamard spaces and in the dynamics of Kleinian groups.
 The main reason for this  place is that  any simply connected,  nonpositively curved space $X$ (a {\em Hadamard space}) has a natural, ``visual'' compactification whose boundary $X(\infty)$ is easily described in terms of {\em asymptotic rays}; and, when $X$ is given a discrete group $G$ of motions, the Busemann functions of rays appear as the densities at infinity of the Patterson-Sullivan measures of $G$ (\cite{roblin}, \cite{sullivan}).

The simple visual picture of the compactification of a Hadamard space unfortunately breaks down for general, non-simply connected manifolds: but Busemann functions (more precisely, their direct generalizations known as {\em horofunctions}) have inspired Gromov to define a natural, universal compactification (the {\em horofunction compactification}), whose properties however are more difficult to describe.
The aim of this paper is to investigate how far the  visual  description of this boundary and the usual properties of rays carry over in the negatively curved, non-simply connected case, and to stress the main differences. 

\vspace{2mm} 
Let us start by describing a first, na\"if approach to the problem of finding a ``good'' geometric compactification of a general complete Riemannian manifold. The first idea is to add all ``asymptotic directions'' to the space, similarly to $\mathbb{E}^n$ which can be compactified as the closed ball $B^n=\mathbb{E}^n \cup S^n$ by adding the set of all oriented half-lines modulo (orientation-preserving) parallelism. Now, on a general Riemannian manifold we have at least two elementary notions of asymptoticity for rays $\alpha, \beta: \mathbb{R}^+ \rightarrow X$ with, respectively,  origins $a,b$:
\vspace{1mm}

\noindent  {\sc \small $\bullet$ Distance Asymptoticity}: we say that  $d_{\infty}(\alpha, \beta) \!<  \! \infty$  if $\sup_t d (\alpha(t), \beta(t) ) \! < \! \infty$, \\ and then we say that $\alpha$ and $\beta$ are  {\em distance-asymptotic} (or, simply,  {\em asymptotic});

\vspace{1mm}
\noindent  {\sc \small  $\bullet$  Visual Asymptoticity}: we say that {\em $\alpha$  tends visually}
\footnote{To avoid an innatural, too  restrictive notion of visual asymptoticity, the correct definition is slightly weaker, cp. \S\ref{corays}, Definition \ref{deficoray}: one allows that $\beta_n = [b_n, \alpha(t_n)]$ for some  $b_n \rightarrow b$. Take for instance a hemispherical cap, with pole $N$, attached to an infinite flat cylinder: two meridians issuing from the pole $N$ (which we obviously want to define the same ``asymptotic direction'')   would never be corays if we do not allow to slightly move the origins of the $\beta_n$.} {\em  to $\beta$},
and write it $\alpha \succ \beta$, if  there exist minimizing geodesic segments $\beta_n=[b,\alpha(t_n)]$ such that
$\beta_n \rightarrow \beta$ (i.e. the angle $\angle \beta, \beta_n \rightarrow 0$);
it is also current to say in this case that  {\em $\beta$ is a coray to $\alpha$}  ($\beta \prec \alpha)$, following Busemann \cite{busemann}.
Then, we say that  {\em $\alpha$ and $\beta$ are visually asymptotic} if $\alpha \succ \beta$ and $\beta \succ \alpha$ ($\alpha \prec \succ \beta$).

\vspace{1mm}
\noindent It is classical that these two notions of asymptoticity coincide for Hadamard spaces. For a Hadamard space  $X$ one   then defines  the {\em visual boundary} $X(\infty)$   as the set of rays ${\cal R} (X)$ modulo asymptoticity, gives to $\overline{X}= X \cup  X(\infty)$ a natural topology  which coincides on $X$ with the original one  and makes of it a compact metrizable space: we will refer to $\overline{X}$ as to the {\em visual compactification} of $X$ (cp. \cite{eberlein} and Section \S\ref{hadamard}).
 \vspace{1mm}

The idea of proceeding analogously for a general Riemannian manifold is tempting but disappointing.
First, besides the case of Hadamard spaces, the relation $\prec$ is known to be generally {\em not symmetric}, and the relation $\prec\succ$ {\em is not an equivalence relation} (exceptly for rays having the same origin, as Theorem \ref{teorcoray} shows).
Some indirect \footnote{The work    \cite{innami}  of Innami concerns the construction of a maximal coray which is not a maximal ray; this property implies that the coray relation is not symmetric.}
examples of the asymmetry can be found in literature for surfaces with {\em variable curvature} \cite{innami}, or for graphs \cite{papa}. We shall give in Section \S\ref{ex} an example of  hyperbolic surface (the Asymmetric Hyperbolic Flute) which makes evident the general asymmetry of the coray relation, which can be interpreted in terms of the geometrical asymmetry of the surface itself. More difficult is to exhibit a case where   $\prec\succ$  is not an equivalence relation: Theorem \ref{teorcoray} and Example \ref{exladder}(a)  (the Hyperbolic Ladder) will make it explicit. 
The problem that visual  asymptoticity is not an equivalence relation has been by-passed by some authors (\cite{lewis}, \cite{nasu}) by taking the  equivalence relation generated by $\prec$ (this means partitioning all the corays to some ray $\alpha$ into maximal packets all of which contain only rays corays to each other): {\em this exactly coincides with taking rays with the same Busemann function} (see \cite{kim} and Section \S2), which explains the original interest of Busemann in this function. 
We will see in \S\ref{corays} that the condition $B_{\alpha} = B_{\beta}$ geometrically simply means  that  {\em we can see $\alpha (t)$ and $\beta(t)$, for $t\gg0$, under a same direction from {\em any} point of the manifold}.
$\hspace{-10mm}$
 \vspace{1mm}

\noindent Secondly, distance and visual asymptoticity (even in this stronger form) are strictly distinct relations on general manifolds: there exist rays staying at bounded distance from each other having different Busemann functions,  and also, more surprisingly, diverging rays defining the same Busemann function. This already happens in constant negative curvature:
 
\begin{theorem}
\label{teordifferentrelations}
{\em 
{\small (The Hyperbolic Ladder \ref{exladder} \& The Symmetric Hyperbolic Flute \ref{exsymmflute})}
}

\noindent There exist hyperbolic surfaces $S_1, S_2$ and rays $\alpha_i, \alpha_i'$ on $S_i$ such that: \\
 (i) $d_{\infty}(\alpha_1, \alpha_1') <  \infty$ but $B_{\alpha_1} \neq B_{\alpha_1'}$; \\
 (ii) $B_{\alpha_2} = B_{\alpha_2'}$ but $d_{\infty}(\alpha_2, \alpha_2') = +  \infty$.
\end{theorem}

\noindent Worst, trying to define a boundary $X^d (\infty)$ or $X^v (\infty)$ from ${\cal R} (X)$ by identifying rays under any of these asymptotic relations generally leads to a non-Hausdorff space, because these relations are not closed (with ${\cal R} (X)$ endowed of the topology of uniform  convergence on compacts):

\begin{theorem}
\label{teornotclosed}
{\em
{\small  (The Twisted Hyperbolic Flute \ref{extwistedflute})}
}

\noindent There exist a hyperbolic surface $X$ and rays $\alpha_n \rightarrow \alpha$ on $X$ such that: \\
 (i) $d_{\infty}(\alpha_n, \alpha_m) <  \infty$ but $d_{\infty}(\alpha_n, \alpha) =  \infty$ for all $n,m$; \\
 (ii) $B_{\alpha_n} = B_{\alpha_m}$ but $B_{\alpha_n} \neq B_{\alpha}$ for all $n,m$.
\end{theorem}
 
\noindent This prevents doing any reasonable measure theory (e.g. Patterson-Sullivan theory) on any compactification built out of    $X^d (\infty)$, $X^v (\infty)$. A remarkable example where this problems occurs is the Teichmuller space ${\cal T}_g$ which, endowed with the Teichmuller metric, has a non-Hausdorff visual boundary for $g \geq 2$  \cite{papa2}.

\vspace{2mm}
Gromov's idea of compactification overrides the difficulty of using asymptotic rays, by considering the topological embedding
\vspace{-2mm}

$$b: X \hookrightarrow C(X) / \mathbb{R} \;\;\; P \mapsto [d(P, \cdot)]$$

\noindent of any Riemannian manifold $X$ in the space of real continuous functions on $X$  (with the  uniform topology), up to additive constants.
He  defines $\overline{X}$ as the closure of $b(X)$ in $C(X)/ \mathbb{R} $, and  its  boundary  as $\partial X = \overline{X} - b(X)$, obtaining a compact, Hausdorff (even metrizable) space where $X$ sits in. 
We will call  $\overline{X}$  the {\em horofunction compactification}
\footnote{This construction first appeared, as far as we know, in \cite{gromov} (cp. also, for instance \cite{bgs}, \cite{bridson}) and for this is also known as the {\em Gromov compactification} (or also as the {\em Busemann} or {\em  metric} compactification) of $X$. We will stick to the name ``horofunction compactification'', keeping the other for the well-known compactification of Gromov-hyperbolic spaces.}
of $X$, and  $\partial X$  the  {\em horoboundary of $X$}. \\
The points of $\partial X$   are  commonly called {\em horofunctions}; 
Busemann functions then naturally arise as particular horofunctions: actually, for points of $X$ diverging along a ray $P_n=\alpha(n)$ we have that 
\vspace{-3mm}

$$b(P_n)=[d(P_n, \cdot)] = [ d(x,P_n) - d(P_n, \cdot)] \longrightarrow [B_{\alpha}(x, \cdot)]$$

\vspace{-1mm}
\noindent in  $C(X) / \mathbb{R}$, cp. Section \S2 for details. 
Accordingly, the  {\em Busemann map}
\vspace{-3mm}

$$B: {\cal R} (X) \rightarrow \partial X$$

\vspace{-1mm}
\noindent  is the map which associates to each ray the class of its Busemann function. For Hadamard manifolds, it is classical that $B$ induces a homeomorphism between the visual boundary  $X(\infty)$ and the horoboundary $\partial X$ (cp. \S2). \\
\vspace{1mm}

The properties of the Busemann map for general nonpositively curved Riemannian manifolds will be the  second object of our interest in this paper. 
The main questions we address are:
\vspace{4mm}

{\em (a) the Busemann Equivalence}: i.e., when do the Busemann functions of two distinct rays  coincide?  
\vspace{1mm}

\noindent Actually, the equivalence relation generated by the coray relation is difficult to test in concrete examples.
In Section \S\ref{busmap} we discuss several notions of equivalence of rays related to the Busemann equivalence; then  we give a  characterization (Theorem \ref{teoreq}) of the Busemann equivalence for rays on  quotients  of Hadamard spaces, in terms of the points at infinity of their lifts, which we call {\em weak $G$-equivalence}. For rays with the same origin, it  can be stated as follows:
\vspace{-2mm}

\begin{criterium}
 \label{critequivalence}
Let  $X=G \backslash \tilde{X}$  be a regular quotient of a Hadamard space. \\
Let  $\alpha, \beta$ be rays based at $o$, with lifts $\widetilde{\alpha}, \widetilde{\beta}$  from $\tilde{o} \in  \tilde{X}$, and let $H_{\tilde{\alpha}}, H_{ \widetilde{\beta}}$ be the horoballs through $\tilde{o}$ centered at the respective points at infinity  $\widetilde{\alpha}^+, \widetilde{\beta}^+$. Then: 
\vspace{-6mm}

$$B_{\alpha} \! = \! B_{\beta} \Leftrightarrow \exists \, (g_n), (h_n) \!\in\! G \mbox{  such that }
\left\{ \! \begin{array}{l}  g_n \widetilde{\alpha}^+ \rightarrow \beta^+ \\ 
		\!\! d( g_n^{-1} \tilde{o} , H_{\tilde{\alpha}} ) \rightarrow 0 \end{array} \right.
\! \! \mbox{and}
 \left\{ \! \begin{array}{l}  h_n \widetilde{\beta}^+ \rightarrow \alpha^+ \\
 		\!\!d( h_n^{-1} \tilde{o}, H_{\tilde{\beta}} ) \rightarrow 0 \end{array} \right.
$$
\end{criterium}

\noindent This reduces the problem of the Busemann equivalence for rays $\alpha, \beta$ on quotients of a Hadamard space  to the problem of approaching the  limit points (of their lifts) $\widetilde{\alpha}^+, \widetilde{\beta}^+$  with sequences $g_n \widetilde{\beta}^+$,  $h_n \widetilde{\alpha}^+$ in the respective orbits, keeping at the same time control of the dynamics of the inverses $g_n^{-1},  h_n^{-1}$.

\vspace{3mm}
{\em (b) the Surjectivity of the Busemann map}: i.e., is any point in the horoboundary of $X$ equal to the Busemann function of some ray?
\vspace{1mm}

\noindent In this perspective, it is natural to extend the Busemann map $B$ to the set  $q{\cal R} (X)$  of {\em quasi-rays} (i.e. half-lines $\alpha : R^+ \rightarrow X$ which are only {\em almost-minimizing}, cp. Definition \ref{defiexcess});  we then call {\em Busemann boundary}  ${\cal B}X = B( q{\cal R} (X) ) $. 
The problem whether ${\cal B}X = \partial X$ has been touched by several authors for surfaces with finitely generated fundamental group (cp.  \cite{shioya} and  \cite{yim}). In \cite{yim} there are examples of a {\em non-negatively} curved surface admitting horofunctions  which are not in ${\cal B}X$, and even of surfaces where the set of Busemann functions of rays emanating from one point is different from that of rays emanating from another \nopagebreak point\footnote{For surfaces with finite total curvature, Yim uses the terminology {\em convex} and {\em weakly convex at infinity}, which is suggestive of the meaning of the value of $2 \pi \chi (X) - \int_X K_K$ \linebreak (to be interpreted, for surfaces with boundary,  as the convexity of the boundary). However this can be misleading, suggesting the possibility of joining with bi-infinite rays any two points at infinity. As our manifolds will generally be infinitely connected, we will not adhere to this terminology.}.
This explains our interest in considering rays with variable initial points, instead of keeping the base point fixed once and for all.

\noindent In \cite{ledrappier} Ledrappier and Wang start developing the Patterson-Sullivan theory on non-simply connected manifolds, and the question whether an orbit  accumulates to a limit point which is a true Busemann function naturally arises; the Theorem below shows that, in this context,  Patterson-Sullivan theory must  take into account limit points which are not Busemann functions, and that some paradoxical facts already happen in the simplest cases:
 \vspace{-2mm}
 
\begin{theorem}
\label{teornonsur}
{\em
{\small (The Hyperbolic Ladder \ref{exladder})}
}

\noindent There exists a Galois covering $X \rightarrow \Sigma_2$ of a hyperbolic surface of genus $2$, with automorphism group  $\Gamma \cong \mathbb{Z}$, such that:

\noindent (i)  ${\cal B}X$ consists of 4 points, while  $\partial X$ consists of a continuum of points;

\noindent (ii) the limit set  $L\Gamma= \overline{\Gamma x_0} \cap \partial X$  depends on the choice of the base point $x_0$, and for some $x_0$ it is included in  $\partial X - {\cal B}X$.
\end{theorem}
\vspace{-1mm}

\noindent The problem of surjectivity and the interest in finding Busemann points in the horoboundary   seems to have been  revitalized due to recent work on Hilbert spaces (cp. \cite{walsh1}, \cite{walsh2}), on the Heisenberg group \cite{klein},  on word-hyperbolic groups and  general Cayley graphs (cp. \cite{bjork},  \cite{webster}). A construction similar to that of Theorem \ref{teornonsur}  is discussed in \cite{bridson} as an example where the {\em boundary of a Gromov-hyperbolic space} does not coincide with the horoboundary (notice however that the notion of boundary for Gromov-hyperbolic spaces  differs from ${\cal B}X$,  as it is defined up to a bounded function).
\vspace{4mm}

{\em (c) the Continuity of the Busemann map}: i.e., how does the Busemann functions change with respect to the initial direction of rays?
\vspace{1mm}

\noindent This is crucial to understand the topology of the horofunction compactification and, beyond the simply connected case, it has not been much investigated in literature so far.
Busemann himself  seemed to exclude it in full generality.\footnote{He wrote in his book  \cite{busemann}: ``{\em It is not possible to make statements about the behaviour of the function $B_{\alpha}$ under general changes of $\alpha$ [...]}''.}  \\
We shall see that, in general,  the dependence from the initial conditions is only lower-semicontinuous:
 \vspace{-1mm}
 
\begin{theorem}
\label{teornoncont}
{\em
{\small (Proposition  \ref{propsemicont} \& The Twisted Hyperbolic Flute \ref{extwistedflute})}
}

\noindent Let  $X=G \backslash \tilde{X}$  be the regular quotient of a Hadamard space.  \\
\noindent (i)  for any sequence of rays $\alpha_n \rightarrow  \alpha$ we have 
$\lim_{n\rightarrow \infty} B_{\alpha_n} \geq B_{\alpha}$; 
 
 \noindent (ii)  there exists  $X= G \backslash \mathbf{H}^2$ and rays $\alpha_n \rightarrow  \alpha$ such that
$\lim_{n\rightarrow \infty} B_{\alpha_n} > B_{\alpha}$.
(Convergence of rays is always meant uniform on compacts.)
\end{theorem}

 \vspace{-1mm}
\noindent The example of the Twisted Hyperbolic Flute \ref{extwistedflute} is the archetype where a jump between $\lim_{n \rightarrow \infty} B_{\alpha_n}$ and $B_{\alpha}$  occurs; 
we will explain geometrically --actually, we produce,  the discontinuity in terms of a  discontinuity in the limit of the {\em maximal horoballs} associated to the $\alpha_n$  in the universal covering, cp. Definition \ref{definonhorospherical} and Remark \ref{remcontinuity}.
 Interpreting $e^{-B_{\alpha}(o, \cdot)}$ as a reparametrized distance to the point at infinity of $\alpha$, 
 the jump can be seen as a  hole suddenly appearing in a limit direction of a hyperbolic sky. \\
 We stress the fact that the problem of continuity makes sense only for {\em rays} $\alpha_n$ (i.e. whose   velocity vectors yield minimizing directions):  it is otherwise easy to produce a discontinuity in the Busemann function  of a sequence of quasi-rays tending to some limit curve which is not minimizing  (and for which the Busemann function may be not defined), cp. Example \ref{exdisc1} in \S\ref{busmap} and the discussion therein.
\vspace{1mm}

It is noticeable that all the possible pathologies in the geometry of rays which we described above already occur for hyperbolic surfaces belonging to  two basic classes: flutes and ladders, see Section  \S\ref{ex}. These are surfaces with infinitely generated fundamental group whose topological realizations are, respectively, infinitely-punctured spheres and $\mathbb{Z}$-coverings of a compact surface of genus $g \geq 2$:
\vspace{-4mm}

\begin{figure}[h]
\centering
\includegraphics[width=80mm, height=30mm]{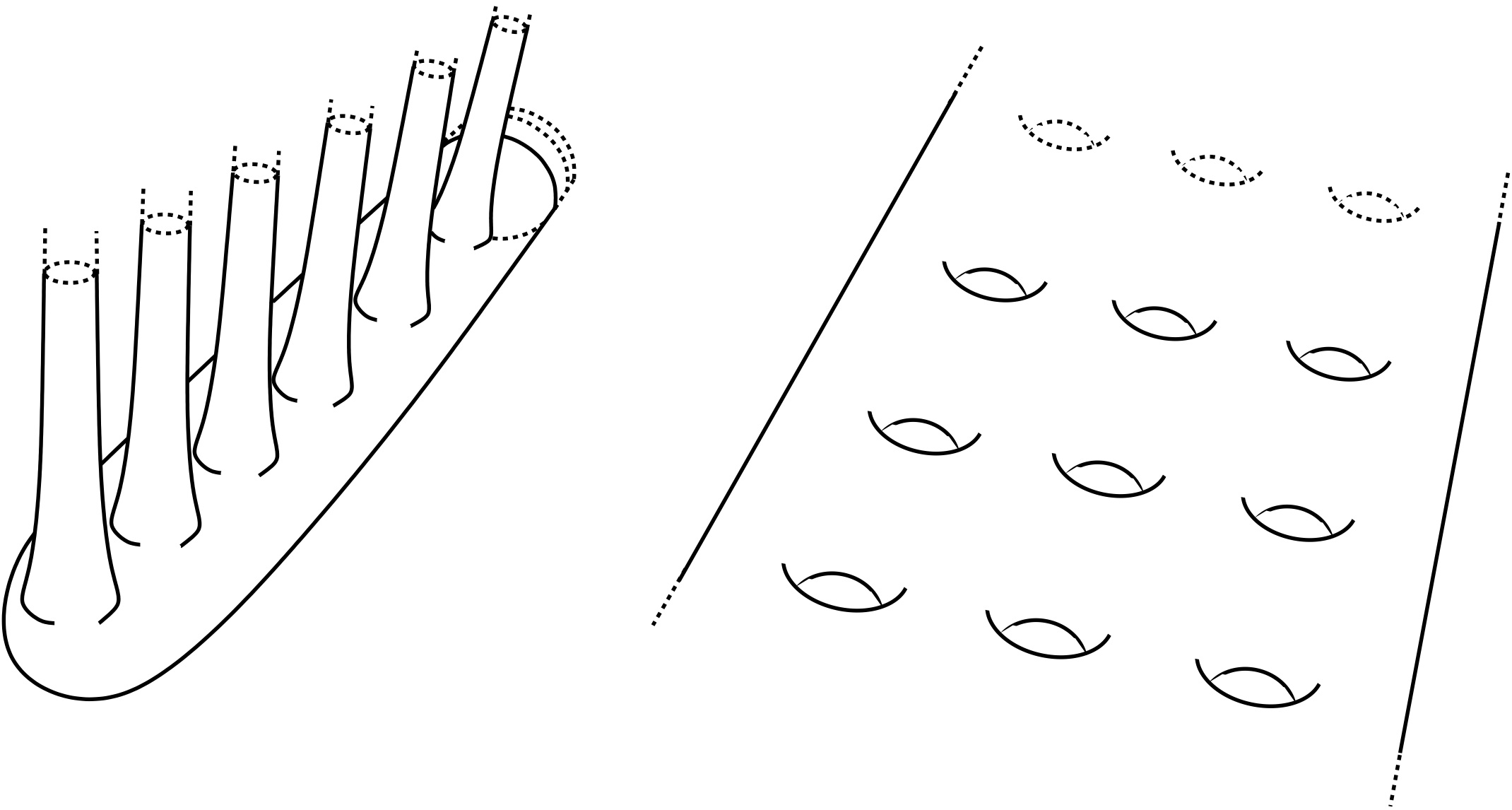}
\caption{{\em Geometric realization of flutes and ladders}}
\label{figflute-ladder}
\end{figure}
 
\vspace{-2mm}
On the other hand, limiting ourselves to the realm of surfaces with finitely generated fundamental group, all the above pathologies disappear and we recover the familiar picture of rays on  Hadarmard manifolds. More generally,  in Section \S\ref{gf} we will consider properties of rays and the Busemann map for {\em geometrically finite manifolds}: these are the geometric generalizations, in dimension greater than $2$, of  the idea of negatively curved surface with finite connectivity (i.e. finite Euler-Poincar\'e characteristic). The precise definition of this class and much of these manifolds is due to Bowditch \cite{bow}; we will summarize the necessary definitions and properties in Section \S\ref{gf}. We will prove: 
\vspace{-1mm}
  
\begin{theorem}
\label{teorgf}
{\em
{\small (Propositions \ref{propprerays}, \ref{propeq}, \ref{propsur}, \ref{propcont} \& Corollary \ref{corgf})}
}

\noindent Let $X= G \backslash \tilde{X}$ be a geometrically finite manifold:

\noindent (i) every quasi-ray on $X$ is finally a ray (i.e. it is a {\em pre-ray}, cp. Definition \nolinebreak \ref{defiexcess});

\noindent (ii) $d_{\infty}(\alpha, \beta) < \infty \; \Leftrightarrow \; B_{\alpha}=B_{\beta} \; \Leftrightarrow \;  \alpha \prec \beta$ for rays $\alpha, \beta$ on $X$;

\noindent (iii) the Busemann map ${\cal R}(X) \rightarrow \partial X$ is surjective and continuous.

\noindent As a consequence, $X(\infty)={\cal R}(X) /_{(Busemann \; eq.)}$ is   homeomorphic to $\partial X$  and:

\noindent $\bullet$ if dim$(X)=2$,    $\overline{X}$  is a compact  surface with boundary;

\noindent $\bullet$  if dim$(X)>2$,  $\overline{X}$  is a compact  manifold with boundary, with a finite number of conical singularities (one for each class of maximal parabolic subgroups of \nolinebreak $G$).
\end{theorem}
 \vspace{-1mm}

\noindent In this regard, it is of interest to recall that  the question whether any geometrically finite manifold has {\em finite topology} (i.e., is homeomorphic to the interior of a compact manifold with boundary) was asked by  Bowditch in \cite{bow}, and recently answered by  Belegradek and Kapovitch, cp. \cite{belekapo}. However, Belegradek and Kapovitch's proof yields a natural topological compactification whose boundary points are  less related to the geometry of the interior than in the horofunction compactification. According to \cite{belekapo},
any horosphere quotient is diffeomorphic to a flat Euclidean vector bundle over a compact base,
so a parabolic end can be seen as the interior of a closed cylinder over a closed disk-bundle.
On the other hand, in the horofunction compactification, a parabolic end is compactified as a cone over the Thom space of this disk-bundle
(cp.Corollary \ref{corgf}\&Example \ref{extomate}); one pays the  geometric content of the horofunction compactification by the appearing of (topological) conical singularities. \\
The problem of relating the ideal boundary and the horoboundary  for geometrically finite groups has also been considered in  
\cite{macpherson}; there, the authors prove  that, in the case of arithmetic lattices of symmetric spaces, both compactifications  coincide with the Tits compactification, and also discuss the relation with the  Martin boundary.

\vspace{3mm}
Section \S\ref{pre} is preliminary: we report here some generalities about the Busemann functions and the coray relation. 

From Section \S\ref{nonpos} on, we focus on nonpositively curved manifolds. 
We briefly recall  the classical visual properties of rays on Hadamard spaces,  and then we turn our attention to their quotients  $X=G \backslash \tilde{X}$. The difference between rays and quasi-rays is deeply related with the different kind of points at infinity of their lift  to $\tilde{X}$; that is why we review a dictionary between limit points of $G$ and corresponding quasi-rays on $X$. Then, we
 prove a formula  ({\bf Theorem \ref{teorbusemann}})  expressing the Busemann function of a ray $\alpha$ on $X$   in terms of the Busemann function of a lift  $\widetilde{\alpha}$ of $\alpha$ to  $\tilde{X}$.
  We will use this formula to translate the Busemann equivalence in terms of  the above mentioned weak $G$-equivalence; this turns out to be the key-tool for constructing examples having Busemann functions with prescribed behaviour.
  
In Section \S\ref{busmap} we discuss the properties of the Busemann map on general quotients of Hadamard spaces; here we prove the Criterium \ref{critequivalence} and the lower semicontinuity.

Section \S\ref{gf} is devoted to  geometrically finite manifolds and to the proof of Theorem \ref{teorgf}.

Finally, we collect in Section \S\ref{ex} the main examples of the paper (the Asymmetric, Symmetric and Twisted Hyperbolic Flute, and the Hyperbolic Ladder).

In the Appendix we report, for the convenience of the reader, proofs of those facts which are either classical, but essential to our arguments, or which we were not able to find easily in literature.
\vspace{3mm}

\begin{minipage}[l]{115mm}
\small \noindent
{\em We will always assume that  geodesics are parametrized by arc-length, and we will use the symbol $[p,q]$ for a  minimizing geodesic segment connecting two points $p,q$. Moreover, we shall often use, in computations,  the notations $x \apprle_{\epsilon} y$   for  $x \leq y +\epsilon$ (respectively,  $x \eqsim_{\epsilon} y$ for $|x-y| \leq \epsilon$)  and abbreviate $d(x,y)$ with  $xy$.}
\end{minipage}

\section{Busemann functions on Riemannian manifolds}
\label{pre}

\subsection{Horofunctions and Busemann functions}
\label{busemannfunctions}

\noindent Let $X$ be any complete Riemannian manifold (not necessarily simply connected). The horofunction compactification of $X$ is obtained  by embedding $X$ in a natural way into the space $C(X)$ of real continuous functions on $X$, endowed with the $C^0$-topology (of  uniform convergence on compact sets):
\vspace{-2mm}

$$ b: X \hookrightarrow C(X) \hspace{10mm} P \mapsto -d(P,\cdot)$$
\noindent  then, defining  $\overline{X} \doteq \overline{b (X)}$ and $\partial X  \doteq \overline{X}- b (X) $. \\
An (apparently) more complicate version of this construction has the advantage of making the Busemann functions naturally  appear  as boundary points. \linebreak
For fixed $P$,  define the {\em horofunction} {\em cocycle} as the function of $x,y$:
\vspace{-2mm}

$$b_{P}(x,y) =  d(x,P) - d(P, y)$$
then,  consider the space  of functions in $C(X)$ up to an additive costant (with the quotient topology) and the same map 
\vspace{-4mm}

$$ b: X \rightarrow C(X)/\mathbb{R} \hspace{10mm} P \mapsto [-d(P,\cdot)]= [ d(x,P) - d(P, \cdot)] = [b_P (x, \cdot) ]$$

\noindent (which is independent from the choice of $x$). The following properties hold, in all generality, for any complete Riemannian manifold, and can be found, for instance, in  \cite{ballmann} or \cite{bridson}:
\vspace{2mm}

\noindent 1) $b$ is a topological embedding,  i.e.  an injective map which is a homeomorphism when restricted to its image;
\vspace{2mm}

\noindent 2) $\overline{X}$ is a compact, 2$^{nd}$-countable, metrizable space.
  
\begin{definition}
\label{defiboundary}
Horoboundary and horofunctions \\
{\em 
The  \emph{horofunction compactification of $X$} and the   {\em horoboundary of $X$}  are respectively the sets $\overline{X} \doteq \overline{b (X)}$ and  $\partial X \doteq \overline{X} - b (X)$.
A \emph{horofunction} is an element $\xi \in\partial X$, that is the limit of a sequence  $[b_{P_n}]$, for $P_n \in X$ going to infinity; we will write $\xi=B_{(P_n)}$.  
}
\end{definition}

\vspace{-1mm}
\noindent Notice that, as $b_{P}(x,y) - b_{P}(x',y) =  b_P (x,x')$  saying that $(P_n) \rightarrow \xi \in \partial X$ is equivalent to  saying that, for any fixed $x$, the  horofunction cocycle $b_{P_n} (x,\cdot)$ converges uniformly on compacts for $n \rightarrow \infty$ (to a representative of $\xi$). \linebreak Concretely, we see a horofunction  $\xi=B_{(P_n)}$ as  a function  of two variables $(x, y)$ satisfying: 
\vspace{2mm}

\noindent (i)  {\em (cocycle condition)} $\;\;B_{(P_n)} (x,y) -  B_{(P_n)} (x',y) = B_{(P_n)} (x,x')$
\vspace{2mm}

\noindent or, equivalently 
\footnote{this formulation is much suggestive as, when thinking of horofunctions as reparametrized distance functions from points at infinity, then we see that the usual triangular inequality becomes an equality  for all points at infinity.}, 
$B_{(P_n)} (x,x') +  B_{(P_n)} (x',y) = B_{(P_n)} (x,y)$. \\
\vspace{1mm}

 The following properties follow right from the definitions:
\vspace{2mm}

\noindent (ii) {\em (skew-symmetry)} \hspace{5mm}  $B_{(P_n)} (x,y) = - B_{(P_n)} (y,x)$
\vspace{2mm}

\noindent (iii) {\em (1-Lipschitz)} \hspace{10mm}  $B_{(P_n)} (x,y) \leq d(x,y)$
\vspace{2mm}

\noindent (iv) {\em (invariance by isometries)}  $\;\;B_{(gP_n)} (gx,gy)=B_{(P_n)} (x,y) \;\; \forall g \in Isom (X)$
\vspace{2mm}

\noindent (v) {\em (continuous extension)} \;  the cocycle $b_P(x,y)$  can be extended to a continuous function $B: X\times \overline{X} \times X \rightarrow \mathbb{R}$, i.e.  $B_{\xi} (x, y) = \lim_{n \rightarrow \infty} b_{P_n} (x,  y)$ if $(P_n) \rightarrow \xi$;
\vspace{2mm}

\noindent (vi) {\em (extension to the boundary)}  every  $g \in Isom (X)$ naturally extends to a homeomorphism $g: \partial X \rightarrow \partial X$.
\vspace{5mm}

  Now, the simplest way of diverging, for  a sequence of points $\{ P_n \}$ on an open manifold $X$, is to go to infinity along a geodesic. As we deal with non-simply connected manifolds, we shall need to distinguish between geodesics and minimizing geodesics:

\begin{definition}
\label{defiexcess}
Excess and quasi-rays

{\em
\noindent The {\em length excess} of a curve $\alpha$ defined on an interval $I$ is the number 
\vspace{-2mm}

$$\Delta(\alpha) = \sup_{t,s\in I} \ell (\alpha; t,s)- d(\alpha(t), \alpha(s))$$
 that is the greatest difference between the length of $\alpha$ between two of its points, and their effective distance.  
Accordingly, we say that a geodesic  $\alpha$ in a manifold $X$ is {\em quasi-minimizing} if $\Delta(\alpha) < +\infty$, and  {\em $\epsilon$-minimizing} if $\Delta(\alpha) \leq \epsilon$. \\ A  {\em quasi-ray}  is  a quasi-minimizing half-geodesic $\alpha: \mathbb{R}_+ \rightarrow X$. 
 For a quasi-ray $\alpha$ there are three possibilities:

\vspace{1mm}
\noindent $\bullet$ either $\alpha$ is minimizing (i.e. $\Delta(\alpha)=0$): then $\alpha$ is a true {\em ray}; 

\vspace{1mm}
\noindent $\bullet$ or $\alpha |_{[t_0, + \infty]}$ is  minimizing for some $a>0$, an then we call $\alpha$ a {\em pre-ray};

\vspace{1mm}
\noindent  $\bullet$ or $\Delta(\alpha) < \infty$ but  $\alpha |_{a, + \infty}$ is never minimizing, for any $a \in \mathbb{R}$; in this case, following \cite{haas},  we call $\alpha$ a {\em rigid quasi-ray}.
\vspace{1mm}

\noindent We will denote by ${\cal R}(X)$ and $q{\cal R}(X)$ the sets of rays and quasi-rays of $X$  (resp. ${\cal R}_o(X)$ and  $q{\cal R}_o(X)$  those with origin $o$), with the uniform topology given by convergence on compact sets.
}
\end{definition}

There exist, in literature, examples of all three kinds of quasi-rays. \\ 
An enlightening  example is the modular surface $X= PSL (2, \Z) \backslash \Hy^2$ (though only an orbifold). $X$ has a 6-sheeted, smooth covering $\hat{X}=\Gamma(2) \backslash \Hy^2 \rightarrow X$, with finite volume; the half-geodesics $\alpha$ of $\hat{X}$ with infinite excess are precisely the bounded geodesics and the unbounded, recurrent ones (those who come back infinitely often in a compact set); their lifts in the half-plane model of $\Hy^2$ correspond to the half-geodesics  $\widetilde{\alpha}$ having extremity $\widetilde{\alpha}^+ \in \R - \Q$ . Moreover, $\alpha$ is bounded if and only if  $\widetilde{\alpha}^+$ is a {\em badly approximated} number (i.e. its continued fraction expansion is a sequence of bounded integers), see \cite{dalbo}. In this case, all half-geodesics $\alpha$ with $\Delta(\alpha)< \infty$ (corresponding to lifts $\widetilde{\alpha}$ with rational extremity) are minimizing after some time, i.e. they are pre-rays. \\
On the other hand, in  \cite{haas} one can find examples and classification of rigid quasi-rays on particular (undistorted) hyperbolic flute surfaces. 
\vspace{3mm}

For future reference, we report here  some properties of the length excess:

\begin{properties}
\label{lemmaexcess}
Let  $\alpha, \alpha_k: [0,+\infty] \rightarrow X$  curves with origins  respectively  $a$, $a_k$:
\vspace{1mm}

\noindent (i) if $\Delta(\alpha) < \infty$,  then for every $\epsilon> 0$ there exists $T_{\epsilon} \gg 0$ such that
\vspace{-3mm}

$$\Delta ( \alpha |_{[T_{\epsilon}, +\infty]} ) \leq \epsilon \;\;\;\;
\mbox{ and } \;\;\;\;
\Delta ( \alpha |_{[0,T_{\epsilon}]} ) \geq \Delta ( \alpha) - \epsilon \; ; $$
\vspace{-4mm}

\noindent (ii) if $\alpha_k \rightarrow \alpha$ uniformly on compacts, then $\Delta(\alpha) \leq \liminf_{k \rightarrow \infty} \Delta( \alpha_k)$.\\
In particular, any limit of  minimizing geodesics segments is minimizing.
\vspace{1mm}

\noindent (iii) Assume now that  the universal covering of $X$ is a Hadamard space. \\
If    $\widetilde{\alpha}$ is a lift of $\alpha$ to  $\tilde{X}$ with origin $\tilde{a}$, then:
\vspace{-1mm}
$$\Delta (\alpha)= \lim_{t \rightarrow +\infty} d(\tilde{a}, \widetilde{\alpha}(t)) -  d(a,  \alpha(t)) $$
\end{properties}

{\em Proof.} (i) follows from the fact that
the excess is increasing with the width of intervals. 
For (ii), pick $T_{\epsilon}$ as in (i) for $\alpha$,  and   $k \gg 0$ such that  $d(\alpha_k (t), \alpha (t)) \leq \epsilon$ for all $t \in [0,T_{\epsilon}]$; then
\vspace{-3mm} 

$$  a_k  \alpha_k (T_{\epsilon})  \apprle_{2\epsilon}  a \alpha (T_{\epsilon})  \apprle_{\epsilon}   T_{\epsilon} -  \Delta (\alpha)   = \ell (\alpha_k) - \Delta (\alpha)  $$ 

\noindent therefore $\Delta (\alpha_k) \geq \Delta (\alpha) - 3 \epsilon$. By passing to limit for $k \rightarrow \infty$, as $\epsilon$ is arbitrary, we deduce $\liminf_{k \rightarrow \infty} \Delta( \alpha_k) \geq \Delta (\alpha) $.
Finally, if $\tilde{X}$ is Hadamard then $d(\tilde{a}, \widetilde{\alpha}(t)) = t =  \ell(\alpha; 0,t) $ for all $t$; hence, by monotonicity of the excess on intervals,  
\vspace{-5mm} 

$$\Delta (\alpha)= \lim_{t \rightarrow +\infty}  \ell(\alpha; 0,t) - d(a, \alpha(t))=  \lim_{t \rightarrow +\infty} d(\tilde{a}, \widetilde{\alpha}(t)) -  d(a,  \alpha(t)) \;  \mbox{ \qed}$$

\begin{proposition}
\label{propbusemannrays}
Let $\alpha: \mathbb{R}_+ \rightarrow X$ be a quasi-ray.
Then,  the horofunction cocycle  $b_{\alpha(t)}(x, y)$ converges uniformly on compacts to a horofunction for $t \rightarrow + \infty$.
\end{proposition}

\begin{definition} Busemann functions\\
{\em  Given a quasi-ray $\alpha$, the cocycle $b_{\alpha(t)}(x, y)$ is called a {\em Busemann cocycle}, and the horofunction $B_{\alpha} (x,y) = \lim_{t \rightarrow \infty} b_{\alpha(t)}(x, y)$ is called a {\em Busemann function}; the Busemann function of $\alpha$ will also be denoted by  $\alpha^+$. \\
The {\em Busemann map} is the map
\vspace{-2mm}

$$B:  q{\cal R}(X) \rightarrow \partial X \hspace{1cm} \alpha \mapsto B_{\alpha}$$

\noindent The image of this map, denoted ${\cal B}X$, is the subset of {\em Busemann functions}, that is those particular horofunctions associated to quasi-rays. We shall  denote  ${\cal B}_o X$  the image  of the Busemann map restricted to to $q{\cal R}_o(X)$.
}
\end{definition}

 The proof of Proposition \ref{propbusemannrays} relies on the
\vspace{-1mm}

\begin{propmono}
\label{propquasimonotone}
\noindent Let $\alpha$ be a quasi-ray from $a$:  for all $\epsilon>0$ there exists $T_{\epsilon}$ such that
$b_{\alpha(s)} (a, y) \apprge_{2\epsilon}  b_{\alpha(t)}(a,   y) \; \forall s>t>T_{\epsilon}$.
\end{propmono}

\vspace{-1mm}
\noindent  Actually, if $\Delta(\alpha)= \Delta$, by the property \ref{lemmaexcess}(i) we have, for $s, t \geq T_{\epsilon}$:
\vspace{-3mm}

$$b_{\alpha(s)} (a, y) -  b_{\alpha(t)}(a,   y) =
[  a \alpha(s)  -  \alpha (s) y  ] -  [  a  \alpha(t) - \alpha (t)  y ]$$

\vspace{-8mm}

$$  \apprge_{2\epsilon}  [ \ell (\alpha|_{[0,s]} )  -  \alpha (s) y ]  -  
[ \ell (\alpha|_{[0,t]})  -  \alpha (t) y ]
\geq \ell(\alpha|_{[t,s]} )  -  \alpha(t)   \alpha(s)  \geq 0.$$

\vspace{-1mm}
\noindent Notice that this is a true monotonicity property when $\alpha$ is a ray.

 \vspace{3mm}
{\em Proof of Proposition \ref{propbusemannrays}.}
As $b_{P}(x,y) - b_{P}(x',y) =  b_P (x,x')$, then the cocycle $b_{\alpha(t)}(x, y)$ converges for $t \rightarrow \infty$ if and only if $b_{\alpha(t)}(x', y)$ converges; we may therefore assume that $x=a$ is the origin of $\alpha$. The  Lipschitz functions $b_{\alpha(t)} (a, \cdot)$ are uniformly bounded on compacts, hence a subsequence $b_{\alpha(t_n)}$ of them converges uniformly on compacts, for $t_n \rightarrow \infty$; but, then, property \ref{propquasimonotone} easily implies that $b_{\alpha(t)}$ must also converge uniformly for $t \rightarrow \infty$ to the same limit, and uniformly.\qed

\subsection{Horospheres and the coray relation}
 \label{corays}
 
If $\xi$ is a horofunction and $x \in X$ is fixed, the sup-level set
\vspace{-3mm}

$$H_{\xi} (x)   =\{ y \; | \; \xi (x,y) \geq 0 \}$$
\vspace{-4mm}

\noindent (resp. the  level set $\partial H_{\xi} (x)  = \{ y \; | \; \xi (x,y) = 0 \}$)
is called the {\em horoball} (resp. the {\em horosphere}) {\em centered at $\xi$}, passing through $x$. \\
If $H_{\xi},H'_{\xi}$ are horoballs centered at $\xi \in \partial X$,  we define the {\em signed distance} to a horoball as
\vspace{-4mm}

$$ \rho (x, H_{\xi}) =  
\left\{ \!\!\!\!  \begin{array}{rl} d(x, \partial H_{\xi}) &\!\!\!\! \mbox{if } x \not\in H_{\xi}(y)\\ 
 				-d(x, \partial H_{\xi}) & \!\!\!\! \mbox{otherwise} \end{array} \right.
\;\;
 \rho (H_{\xi}, H'_{\xi}) =  
\left\{ \!\!\!\!   \begin{array}{rl}   d(\partial H_{\xi}, \partial H'_{\xi}) &\!\!\!\! \mbox{if } H_{\xi} \supset H'_{\xi}\\ 
 			          -d(\partial H_{\xi}, \partial H'_{\xi}) &\!\!\!\! \mbox{otherwise}   \end{array} \right.$$
 
\vspace{-1mm}
\noindent By the Lipschitz condition, we always have $B_{\xi} (x,y) \leq \rho  (  H_{\xi} (x),   H_{\xi} (y) )$.

\noindent On the other hand, notice that when $\alpha$ is a ray and  $x=\alpha(t),y=\alpha(s)$ are points on  $\alpha$ with $s>t$ we have
\vspace{-3mm}

$$ B_{\alpha} (x,y) = d ( x,   y )=  \rho (  H_{\alpha^+} (x),   H_{\alpha^+} (y) )$$

\vspace{-1mm}
\noindent It is a remarkable rigidity property that the equality holds precisely for points which lie on rays, which are {\em  corays to $\alpha$}:  
 
\begin{definition}{Corays}
\label{deficoray}

{\em
\noindent The definition of coray  formalizes the idea of  seeing (asymptotically) two rays under the same direction, from the  origin of one of them. \\
A half-geodesic $\alpha$  with origin $a$ is a {\em coray}
\footnote{We stress the fact that, by Property \ref{lemmaexcess}(ii) of the excess, every coray is necessarily a ray.}
to a quasi-ray $\beta$ in $X$ -- or, equivalently, {\em $\beta$  tends visually to $\alpha$} --  (in symbols: $\alpha \prec \beta$)  if there exists a sequence of minimizing geodesic segments $\alpha_n= [a_n,b_n]$ with $a_n \rightarrow a$ and  $b_n=\beta(t_n) \rightarrow \infty$ such that $\alpha_n \rightarrow \alpha$ uniformly on compacts; equivalently, such that $\alpha_n'(0) \rightarrow \alpha'(0)$. \\
If $\alpha \prec \beta$ and $\beta \prec \alpha$,  we write $\alpha \prec\succ \beta$ and say that they are {\em visually asymptotic}.\\
We will say that {\em $\alpha, \beta$ are visually equivalent from $o$} if there exists a ray $\gamma$ with origin $o$ such that $\gamma \prec \alpha$ and $\gamma \prec \beta$ (i.e. if we can see $\alpha$ and $\beta$ under a same direction from $o$).
}
\end{definition}

\vspace{-2mm}

\noindent Given   $x, y \in X$, we denote by  $\overrightarrow{xy}$ a complete half-geodesic which is the continuation,  beyond $y$,  of a {\em minimizing} geodesic segment $[x,y]$. Then:

\begin{proposition}
\label{propcoray}
For any quasi-ray  $\beta$ we have: 
 $\;B_{\beta} (x, y) = d(x, y)\; \Leftrightarrow \; \overrightarrow{xy}\prec \beta$. \linebreak
 In particular, if $B_{\beta} (x, y) = d(x, y)$, the extension of any minimizing segment $[x,y]$ beyond $y$ is always a ray.
 \end{proposition}
 
\begin{remarks}
\label{remcorcoray}
{\em It follows that:

 \noindent (i)  {\em any coray  $\alpha \prec \beta$} (and $\beta$ itself, if it is a ray) {\em minimizes the distance between the $\beta$-horospheres  that it meets};
 
 \noindent (ii)  {\em for any quasi-ray  $\beta$, we have the equality $B_{\beta} (x,y) = \rho  (  H_{\beta} (x),   H_{\beta} (y) )$} \linebreak
 (as it is always possible to define a coray $\alpha$ to $\beta$ intersecting  $H_{\beta} (x)$ and  $H_{\beta} (y)$, and  $B_{\beta}$  increases exactly as $t$ along $\alpha(t)$).
}
\end{remarks}

\begin{theorem}
\label{teorcoray}
\hspace{1mm}  Assume that  $\alpha, \beta$ are rays in $X$ with origins $a,b$ respectively.  \linebreak
 The following conditions are equivalent: 
\vspace{1mm}

\noindent (a) $B_{\alpha} (x,y) \! = \! B_{\beta} (x,y)$  $ \forall x, y \in X$;
\vspace{1mm}

\noindent (b)  $\alpha \prec\succ \beta$ and $B_{\alpha} (a,b) = B_{\beta} (a,b)$;
\vspace{1mm}

\noindent (c) $\alpha, \beta$ are visually equivalent from every $o \in X$.
\end{theorem}

\noindent  Proposition \ref{propcoray} is folklore (under the unnecessary, extra-assumption that $\overrightarrow{xy}$ {\em is} a ray), and it is already present in Busemann's book \cite{busemann}. Theorem \ref{teorcoray} (a)$\Leftrightarrow$(c) is a reformulation in terms of visibility of the equivalence, proved in \cite{kim},  between Busemann equivalence and the coray relation generated by $\prec$;  the part (a)$\Leftrightarrow$(b) stems from the work of Busemann \cite{busemann} and Shiohama \cite{shiohama}, but we were not able to find it explicitly stated anywhere. For these reasons, we report the proofs of both results  in the Appendix.

\begin{remarks}
\label{remcoray}
{\em
${}$
 
\noindent (i) {\em The coray relation is not symmetric} and {\em  the visual asymptoticity is not transitive, in general}, already for (non simply-connected) negatively curved surfaces, as  we will see in the Examples \ref{exasymmflute} and \ref{exladder}. On the other hand, visual asymptoticity is an equivalence relation when restricted to rays having all the same origin, by the above theorem.
\vspace{1mm}

\noindent (ii)  {\em The condition $B_{\alpha} (a,b) = B_{\beta} (a,b)$ is not just a normalization condition}. \\ In Example \ref{exladder} we will show that  there exist rays $\alpha, \beta$ satisfying $\alpha \prec\succ \beta$ , but  such that $B_{\alpha} \neq B_{\beta}$ an do not differ by a constant. 
\vspace{1mm}

\noindent (iii)  Horospheres are generally not smooth, as  Busemann functions and horofunctions generally are only Lipschitz  (cp. \cite{eberlein},\cite{yim})
This explains the possible existence of multiple corays, from one fixed point,  to a given ray $\alpha$, as well as the asymmetry of the coray relation; actually, in every point of differentiability of $B_{\alpha}$, the direction of a coray to $\alpha$ necessarily coincides with the gradient of $B_{\alpha}$ by Proposition \ref{propcoray}.
}
\end{remarks}

 \section{Busemann functions in nonpositive curvature}
 \label{nonpos}

 \subsection{Hadamard spaces}
 \label{hadamard}

Let $\tilde{X}$ be a simply connected, nonpositively curved manifold (i.e. a {\em Hadamard space}).
In this case, every geodesic  is  minimizing; moreover, as the equation of geodesics has solutions which depend continuously on the initial conditions,  ${\cal R} ( \widetilde{X})$ can be topologically identified with the unit tangent bundle $S \tilde{X}$.
\vspace{-1mm}

\begin{proposition}
\label{prophadamard}
Let  $\widetilde{X}$ be  a Hadamard space:
\vspace{1mm}

\noindent (i) if $\alpha, \beta$ are rays, then $d_{\infty} (\alpha, \beta) < \infty \;  \Leftrightarrow \; B_{\alpha}=B_{\beta} \; \;  \Leftrightarrow \; \alpha \prec \beta$.

\noindent  Moreover, two rays with the same origin are Busemann equivalent iff they coincide, so the restriction of the Busemann map 
$B_o \!:\! {\cal R}_o  ( \widetilde{X}) \rightarrow \partial \widetilde{X}$ is injective; \\
{\em (accordingly, we will denote by $[o,\xi]$ the only geodesic starting at $o$ with point at infinity \nolinebreak $\xi$)}
\vspace{1mm}

\noindent (ii) for any $o \in \tilde{X}$, the restriction   of the Busemann map $B_o: {\cal R}_o ( \widetilde{X})\rightarrow \partial \widetilde{X}$ is surjective, hence ${\cal B}\tilde{X} = {\cal B}_o \tilde{X} = \partial \tilde{X}$;
\vspace{1mm}

\noindent (iii) the Busemann map $B: {\cal R}(\tilde{X}) \rightarrow \partial \tilde{X}$ is continuous.
\vspace{1mm}

\noindent The space ${\cal R}_o( \widetilde{X}) \cong S_o(\widetilde{X})$ being compact, the map $B_o$ gives a homeomorphism $S_o(\widetilde{X}) \cong \partial \tilde{X}$ for any $o$
{\em  (for this reason the topology of the horoboundary $\partial \widetilde{X}$ for Hadamard manifolds is  also known as the   {\em sphere} topology). 

\noindent Also notice that  Proposition \ref{propcoray}(a), together with point (a) above, imply the following fact (which we will frequently use):} 
\vspace{1mm}

\noindent (iv) if $B_{\beta} (x,y) = d(x,y)$ for some $x\neq y$, then  $\overrightarrow{xy}^+=\beta^+$.
\end{proposition}

 \noindent The above properties of rays on a Hadamard space are well-known (cp. \cite{ballmann},  \cite{eberlein}, \cite{bridson}); we shall give in the Appendix a unified proof of (i), (ii) and (iii) for the convenience of the  reader.
 Here we just want to stress that the distinctive feature of a Hadamard space which makes this case so special: for any ray $\alpha$,   {\em the Busemann function  $B_{\alpha} (x,y)$ is uniformly approximated on compacts by its Busemann cocycle $b_{\alpha(t)}(x,y)$}. Namely: 
  \vspace{-1mm}
  
\begin{approlemma} 
\label{lemmaapprox}
Let $\widetilde{X}$ be a Hadamard space. \\
For any compact set $K$ and $\epsilon>0$, there exists a function $T (K, \epsilon)$ such that for any $x,y \in K$ and any ray $\alpha$ issuing from $K$,  we have  $|B_{\alpha} (x,y) - b_{\alpha(t)}(x,y)| \leq \epsilon$, provided that $t \geq T (K, \epsilon)$.
\end{approlemma}
\vspace{-1mm}

\noindent In fact, properties  (ii) and (iii)   follow directly from the above approximation lemma, while  (i) is a consequence of  convexity of the distance function on a Hadamard manifold and of standard comparison theorems (cp. \S\ref{apphadamard} for details).
 
\noindent  {\em A uniform approximation result  as Lemma \ref{lemmaapprox} above does not  hold for general quotients of Hadamard spaces}: actually, from a uniform approximation of the Busemann functions by the Buseman cocycles one easily deduces surjectivity and continuity of the Busemann map  as in the proof of (ii)\&(iii) in \S\ref{apphadamard}, whereas Example \ref{exladder} shows that for general quotients of Hadamard spaces  the Busemann map is not surjective.

\subsection{Quotients of Hadamard spaces}
\label{quotients}

Let $X= G \backslash \tilde{X}$ be a nonpositively curved manifold, i.e. the quotient of a Hadamard space by a discrete, torsionless group of isometries $G$ (we call it a  {\em regular} quotient). In this section we explain the relation between the Busemann function of a quasi-ray $\alpha$ of $X$ and the Busemann function of a lift $\widetilde{\alpha}$ of $\alpha$ to $\tilde{X}$, which will be crucial for the following sections. \\
 Let us recall some terminology:

\begin{definition}
\label{definonhorospherical}
{\em Let $G$ be a discrete group of isometries  of a Hadamard space \nolinebreak $\tilde{X}$. 
The {\em limit set of $G$} is the set $LG$ of accumulation points in  $\partial \tilde{X}$ of any orbit  $G \tilde{x}$ of $G$; the set  $Ord \, G = \partial \tilde{X} - LG$ is  the {\em discontinuity domain} for the action of $G$ on $\partial \tilde{X}$, and its points are called {\em ordinary points}.
A point $\xi  \in LG$ is called:

\vspace{1mm}
\noindent $\bullet$ a {\em radial}  point if one (hence, every) orbit $G\tilde{x}$ meets  infinitely many times  an $r$-neighbourhood  of $[x,\xi]$   (for some $r$ depending on $\tilde{x}$); 

\vspace{1mm}
\noindent $\bullet$  a {\em horospherical} point if one (hence,  every) orbit $G\tilde{x}$ meets every horoball  centered at $\xi$,  i.e.   $\sup_{g \in G} B_{\xi} (\tilde{x},g\tilde{x}) = +\infty$ for every $\tilde{x} \in \tilde{X}$.
 }
\end{definition}

\vspace{-1mm}
\noindent  Radial points clearly are horospherical points, and correspond to the extremities of rays $\widetilde{\alpha}$ whose projections $\alpha$ to $X$  come   back infinitely many times into some compact set (so $\Delta(\alpha)=\infty$).
A simple example of non-horospherical point is the fixed point of a parabolic isometry of a Fuchsian group
\footnote{On the other hand, in dimension $n \geq 3$ parabolic points can be horospherical, cp. \cite{starkov}.}
(a {\em parabolic} point).
For finitely generated Fuchsian groups, it is known that all horospherical points are radial, but  starting from dimension $3$ there exist examples of horospherical non-radial  (even parabolic) points (see \cite{dalbo}, \cite{dalbostarkov}).

\vspace{2mm}
If $\xi$ is  non-horospherical, then for every $\tilde{x}$ there exists a {\em maximal horoball} 
\vspace{-2mm}

$$H_{\xi}^{max} (x) = \{ \tilde{y} \in \tilde{X} \; |\; B_{\xi} (\tilde{x}, \tilde{y}) \geq   \sup_{g \in G} B_{\xi} (\tilde{x},g\tilde{x})  \}$$

\vspace{-2mm}
\noindent (only depending on $\xi$ and on the projection $x$ of $\tilde{x}$ on $X= G \backslash \tilde{X}$) whose interior does not contain any point of $G\tilde{x}$.
For Kleinian groups, there is   large freedom in the orbital approach of the maximal horosphere, which leads to the following distinction:

\begin{definition}
{\em Let $\xi$ be a non-horospherical point of $G$, and $\tilde{x} \in \tilde{X}$:

\vspace{1mm}
\noindent $\bullet$  $\xi$ is a {\em $\tilde{x}$-Dirichlet point} if $\tilde{x} \in  H_{\xi}^{max} (x)$,   i.e. $\sup_{g \in G} B_{\xi} (\tilde{x},g\tilde{x}) = 0$;  

\vspace{1mm}
\noindent $\bullet$  $\xi$ is a {\em $\tilde{x}$-Garnett point} if it   is not  $g\tilde{x}$-Dirichlet for all $g\in G$, which means that  \linebreak   $B_{\xi} (\tilde{x},g \tilde{x}) < sup_{g \in G} B_{\xi} (\tilde{x},g\tilde{x}) < +\infty$ for all $g \in G$;

\vspace{1mm}
\noindent $\bullet$ $\xi$ is  {\em universal Dirichlet} if  $\forall \tilde{x} \in  \tilde{X}$   $\exists g \in G$ such that $\xi$ is $g\tilde{x}$-Dirichlet, and a  {\em Garnett point} otherwise.
}
\end{definition}

\noindent In literature  one can find examples of limit points which are $\tilde{x}$-Dirichlet points but $\tilde{x}'$-Garnett for $\tilde{x}' \neq \tilde{x}$, and also of points which are $\tilde{x}$-Garnett for all $\tilde{x}$,  cp. \cite{nicholls1}, \cite{nichollswaterman1}.
Notice that Dirichlet points may be ordinary or limit points; on the other hand, any ordinary point is universal Dirichlet (as if there exists a sequence $g_n \in G$ such that $d(g_n \tilde{x},H_{\xi}^{max} (x)) \rightarrow 0$, then $\xi$ is necessarily a limit point). We will meet another relevant class of universal Dirichlet points in Section \S\ref{gf} (the {\em bounded parabolic points}). 
Notice that we have, by definition:
\vspace{-3mm}

$$ LG = L^{hor}G \, \sqcup \, L^{u.dir}G  \,\sqcup \, L^{gar}G $$
  
\vspace{-2mm}
\noindent a disjoint union of the subsets of horospherical, universal Dirichlet and Garnett \nolinebreak points.
\vspace{2mm}
  
Consider now the  closed {\em Dirichlet domain of $G$  centered at $\tilde{x} \in \tilde{X}$}:
\vspace{-3mm}

$$D(G,\tilde{x}) =\{ y \in \tilde{X} \; | \; d(y,x) \leq d(y, g \tilde{x}) \; \forall g \in G \}$$ 

\vspace{-1mm}
\noindent This is a convex, locally finite
\footnote{i.e. for any compact set $K \subset \tilde{X}$ one has  $gD(G,\tilde{x}) \cap K \neq\emptyset$ only for finitely many $g \in G$.}
fundamental domain for the $G$-action on $\tilde{X}$; we will denote by
\vspace{-5mm}

$$\partial D(G,\tilde{x})= \overline{D(G,\tilde{x})} \cap \partial \tilde{X}$$ 

\vspace{-1mm} 
\noindent its trace at infinity. Then, we have the following characterization, which explains the name   ``Dirichlet point'':

\begin{proposition}{(Characterization of  Dirichlet points)}
{\em
\hspace{1mm}Let $\xi \in \partial \tilde{X}$ and $\tilde{x} \in \tilde{X}$. \\
Then, \emph{ $\xi$ is $\tilde{x}$-Dirichlet if and only if $\xi$ belongs to $\partial D(G,x)$}. 
}
\end{proposition}

{\em Proof.} Let $\tilde{\gamma}= [\tilde{x},\xi]$. As the Dirichlet domain is convex, we have that  $\xi \in \partial D(G, \tilde{x})$ if and only if $\tilde{\gamma}(t) \in D(G, \tilde{x})$ for all $t$, which means that
\vspace{-3mm}

\begin{equation}
\label{conditiondistance}
d(\tilde{\gamma}(t), \tilde{x}) \leq d(\tilde{\gamma}(t), g \tilde{x}) \mbox{ for } t\geq 0  \mbox{ and for all }  g \in G
\end{equation}

\noindent On the other hand, condition (\ref{conditiondistance}) is equivalent to 
\vspace{-3mm}

\begin{equation}
\label{conditionbusemann}
\sup_{g \in G}B_{\xi} (\tilde{x}, g \tilde{x}) \leq 0  \mbox{ (i.e. $\xi$ is $\tilde{x}$-Dirichlet)}
\end{equation}

\vspace{-1mm}
\noindent In fact, we obtain  (\ref{conditionbusemann}) from (\ref{conditiondistance}) by passing to limit for $t \rightarrow +\infty$. Conversely,  (\ref{conditionbusemann}) implies that  $\tilde{x} \in H^{max}_{\xi} (x)$, and as we know that the direction $\tilde{\gamma}$  is the shortest to travel out of the horoball from $\tilde{\gamma}(t)$, we deduce   (\ref{conditiondistance}).\qed
\vspace{3mm}

The relation with the excess is explained by the following:

\begin{excesslemma}
\label{lemmadelta}
Let $X\!\!= \!G \backslash\! \tilde{X}$ be a regular quotient of a Hadamard space  \nolinebreak $\tilde{X}$. \\
Assume that  $\alpha$ is a half-geodesic in $X$ with origin $a$,  and lift it to $\widetilde{\alpha} $ in $\tilde{X}$  with origin $\tilde{a}$. Then:
\vspace{-5mm}

$$ \Delta (\alpha) = \sup_{g \in G} B_{\alpha} (\tilde{a}, g\tilde{a}) = d(\tilde{a}, H^{max}_{\tilde{\alpha}^+} (a))$$
\end{excesslemma}

\emph{Proof.}  
We have, for any $g \in G$:
\vspace{-6mm}

$$\Delta (\alpha) =\lim_{t \rightarrow \infty}  d(\tilde{a}, \widetilde{\alpha}(t)) - d (a, \alpha(t)) \geq 
\lim_{t \rightarrow \infty}  d(\tilde{a},  \widetilde{\alpha}(t)) - d (\widetilde{\alpha}(t), g \tilde{a}) = B_{\tilde{\alpha}} (\tilde{a}, g\tilde{a})$$

\vspace{-2mm}
\noindent so $\Delta(\alpha) \geq \sup_{g \in G} B_{\tilde{\alpha}} (\tilde{a}, g \tilde{a})$.
On the other hand, for arbitrary $\epsilon>0$, let  $t \gg 0$ such that
$\Delta(\alpha|_{[0,t]}) \eqsim_{\epsilon} \Delta (\alpha)$, and let  $g_t\in G$ such that 
$d(a,\alpha(t)) = d (g_t\tilde{a}, \widetilde{\alpha}(t))$.
Then, by monotonicity of the Busemann cocycle  (\ref{propquasimonotone}), we have  for $s> t$ 
\vspace{-2mm}

$$ \Delta (\alpha) \eqsim_{\epsilon}   d(\tilde{a}, \widetilde{\alpha}(t)) - d ( \widetilde{\alpha}(t), g_t\tilde{a})
\leq  d(\tilde{a}, \widetilde{\alpha}(s)) - d ( \widetilde{\alpha}(s), g_t \tilde{a}).$$

\noindent Letting $s\rightarrow +\infty$ we get $\Delta(\alpha) \apprle_{\epsilon} B_{\tilde{\alpha}} (\tilde{a}, g_t \tilde{a})$ and, as $\epsilon$ is arbitrary, we deduce the converse inequality
$\Delta(\alpha) \leq \sup_{g \in G} B_{\tilde{\alpha}} (\tilde{a}, g\tilde{a})$. \\
To show that  $\sup_{g \in G} B_{\alpha} (\tilde{a}, g\tilde{a}) = d(\tilde{a}, H^{max}_{\tilde{\alpha}^+} (a))$ we just notice 
that,   if $\tilde{y}$ is the point cut on $[a,\widetilde{\alpha}^+]$ by $H^{max}_{\tilde{\alpha}^+} (a)$ then, by Proposition \ref{propcoray},
$$d(\tilde{a}, H^{max}_{\tilde{\alpha}^+} (a))= d(\tilde{a}, \tilde{y}) = B_{\tilde{\alpha}} (\tilde{a}, \tilde{y}) 
= \sup_{g \in G} B_{\tilde{\alpha}} (\tilde{a}, g\tilde{a})  \;.\; \mbox{\qed}$$

\begin{theorem}
\label{teorbusemann}
Let $X = G \backslash \widetilde{X}$ be a regular quotient of a Hadamard space $\tilde{X}$. \\
Assume that   $\alpha$ is a quasi-ray on $X$ with origin $a$, and lift it to $\widetilde{\alpha} $ in $\tilde{X}$, with origin $\tilde{a}$. Then, for all $x,y \in X$ we have:
\vspace{-2mm}

\begin{equation}
\label{eqbusemann}
 B_{\alpha} (x,y) = \rho (H_{\tilde{\alpha}^+} ^{max}(x),  H_{\tilde{\alpha}^+} ^{max}(y))
\end{equation}

\noindent In the particular case where $x=a$ the formula becomes:
\vspace{-2mm}

\begin{equation}
\label{eqbusemanna}
 B_{\alpha} (a,y) =  \sup_{g \in G} B_{\tilde{\alpha}} (\tilde{a}, g \tilde{y}) - \Delta ( \alpha )
= \rho ( \tilde{a},   H_{\tilde{\alpha}^+} ^{max}(y)) -  \Delta (\alpha)
\end{equation}

\vspace{-2mm}

\begin{equation}
\label{eqbusemannb}
 \hspace{-25mm} \mbox{ and, if  $\alpha$ is a ray: \hspace{5mm}} B_{\alpha} (a,y) = \sup_{g \in G} B_{\tilde{\alpha}} (\tilde{a}, g \tilde{y}) 
=  \rho ( \tilde{a},   H_{\tilde{\alpha}^+} ^{max}(y))
\end{equation}

\end{theorem}

\noindent 
Notice that, in the particular case of a ray $\alpha$,  formula (\ref{eqbusemanna}) is quite natural, if we interpret $B_{\alpha} (a, y)$ as a (renormalized, sign-opposite) ``distance to the point at infinity'' $\alpha^{+}$ in  $\partial X$; in fact,  the distance on the quotient manifold $X=G \backslash \tilde{X}$ can always be expressed as
$ d(a, y) =    \inf_{g \in G} d(\widetilde{a} , g \widetilde{y}) $.
\vspace{4mm}

\emph{Proof of Theorem \ref{teorbusemann}.}  
We shall first prove the particular formula (\ref{eqbusemanna}). \\
Since $\ell ( \alpha; 0,t) - d (a, \alpha (t)) \leq \Delta (\alpha)$ for all $t$, we get
\vspace{-6mm}

\[ B_{\alpha} (a,y) = \lim_{t} [d(a, \alpha (t)) - d (\alpha (t), y)] \geq \lim_t [\ell (\alpha; 0,t) - \Delta (\alpha) - 
\inf_{g \in G} d (\widetilde{\alpha} (t), g\tilde{y}) ] \]
\vspace{-5mm}
\[  \geq  \lim_t [ d (\tilde{a}, \widetilde{\alpha} (t)) -   d (\widetilde{\alpha} (t), g \tilde{y})] -  \Delta (\alpha)    = B_{\tilde{\alpha}} (\tilde{a}, g \tilde{y}) -  \Delta (\alpha)  \hspace{3mm}\]

\vspace{-1mm}
\noindent for all $g \in G$.  To prove the converse inequality, pick  for each $t>0$ a preimage $\tilde{y}_t$  of $y$ in $\tilde{X}$ such that $d(\alpha(t), y)= d (\widetilde{\alpha}(t), \tilde{y}_t)$.  By monotonicity and Lemma \ref{lemmaexcess}(i) we have, for all $s>t \gg 0$ 
\vspace{-5mm}

$$d(\tilde{a}, \widetilde{\alpha} (s)) - d ( \widetilde{\alpha} (s), \tilde{y}_t) \geq
d(\tilde{a}, \widetilde{\alpha} (t)) - d ( \widetilde{\alpha} (t), \tilde{y}_t) \apprge_{\epsilon}
d(a, \alpha(t)) + \Delta (\alpha) - d(\alpha(t),y)$$

\vspace{-2mm}
\noindent Therefore letting $s \rightarrow + \infty$ we get
\vspace{-3mm}

$$ \sup_{g \in G} B_{\tilde{\alpha}} ( \tilde{a}, g \tilde{y}) \geq  B_{\alpha} ( \tilde{a},\tilde{y}_t) \apprge_{\epsilon}
b_{\alpha(t)} (a, y)+ \Delta(\alpha)$$

\vspace{-1mm}
\noindent and as $\epsilon$ is arbitrarily small, for $t \rightarrow +\infty$ this yields  (\ref{eqbusemanna}). 
Then,  (\ref{eqbusemann}) follows from  (\ref{eqbusemanna}) and the cocycle condition  since, for any   $\tilde{x}' \in \partial H_{\tilde{\alpha}^+} ^{max}(x)$ and $ \tilde{y}' \in \partial H_{\tilde{\alpha}^+} ^{max}(y)$
\vspace{-5mm}

$$B_{\alpha} (x,y) = B_{\alpha} (a,y)- B_{\alpha} (a,x) =B_{\alpha} (\tilde{a},\tilde{y}')- B_{\alpha} (\tilde{a},\tilde{x}') = B_{\tilde{\alpha}} (\tilde{x}', \tilde{y}')  $$

\vspace{-2mm}
\noindent and this  is precisely   the signed distance between the two maximal horospheres, by Remark \ref{remcorcoray}(ii). The second inequalities in  (\ref{eqbusemanna})\&(\ref{eqbusemannb}) are just   geometric reformulations, as for $\tilde{y}' \in H^{max}_{\tilde{\alpha}^+} (y)$ we have
$\sup_{g \in G} B_{\tilde{\alpha}} (\tilde{a}, g \tilde{y}) = B_{\tilde{\alpha}} (\tilde{a}, \tilde{y}') =
d(\tilde{a}, H^{max}_{\tilde{\alpha}^+} (y))$.\qed

\vspace{3mm}
We conclude the section  mentioning the relation between boundary points and type of quasi-rays, which is an immediate Corollary of the Excess Lemma \ref{lemmaexcess}; this was first pointed out by Haas in \cite{haas}, for Kleinian groups of $\Hy^n$:

\begin{corollary}
\label{cordictionary}
 Let $\pi: \tilde{X} \rightarrow X= G \backslash \tilde{X}$ and $\xi \in \partial \tilde{X}$.
  \vspace{1mm}
 
 \noindent (i)  $\xi$ is non-horospherical iff   $\forall \tilde{x}\in \tilde{X}$ the projection $\pi ([ \tilde{x}, \xi])$ is a quasi-ray;
  \vspace{1mm}
 
 \noindent (ii)   $\xi$ is  $\tilde{x}$-Dirichlet iff $\pi ([ \tilde{x}, \xi])$  is a ray; 
 \vspace{1mm}
  
\noindent (iii)   $\xi$ is  $\tilde{x}$-Garnett iff   $\forall g \in G$  the curve  $\pi ([ g\tilde{x}, \xi])$  is a quasi-ray but not a ray.
\end{corollary}
 
\noindent We shall see in Section \S\ref{gf} a special class of manifolds where every quasi-ray is a pre-ray¤: the {\em geometrically finite} manifolds.

\section{The Busemann map}
\label{busmap}

\subsection{The Busemann equivalence}
\label{eq}

 We will consider several different types of equivalence between rays and quasi-rays on quotients of Hadamard spaces. The main motivation for this is to find workable criteria to know when two rays $\alpha, \beta$ are {\em Busemann equivalent}, that is when $B_{\alpha}= B_{\beta}$. 
We first consider the most natural notion  of aymptoticity:

\begin{definition}{Distance asymptoticity}\\
{\em  For quasi-rays $\alpha, \beta$ on a general manifold $X$ we define
\vspace{-2mm}

  $$d_{\infty} (\alpha, \beta) =   \frac12 \limsup_{t \rightarrow +\infty}  \left[ d(\alpha(t), \beta) +  d(\alpha, \beta(t))\right]$$ 

\vspace{-1mm}
\noindent and we say that   $\alpha, \beta$ are \emph{asymptotic} if $d_{\infty}(\alpha, \beta) < \infty$ (resp.  \emph{strongly asymptotic} if $d_{\infty}(\alpha, \beta)=0$); we say that $\alpha, \beta$ are {\em diverging}, otherwise.
}
\end{definition}

\noindent  Notice that  {\em strongly asymptotic quasi-rays define the same Busemann  function}, since for all $\epsilon>0$ there exists $t,s \gg 0$ such that $| b_{\alpha(t)}(x,y) -  b_{\beta(s)}(x,y) | < \epsilon$. \\
On Hadamard spaces  we know, by Proposition \ref{prophadamard}(a), that two rays are Busemann equivalent precisely when they are asymptotic (moreover, for Hadamard spaces of strictly negative curvature, the notions of  asymptoticity and strong asymptoticity coincide).  
Unfortunately, this easy  picture is false in general: the Example \ref{exladder} in Section \S\ref{ex} exhibits, in particular, {\em two asymptotic rays on a hyperbolic surface yielding different Busemann functions}; on the other hand, in the Example \ref{exsymmflute}  we produce a hyperbolic surface with  {\em two diverging rays defining the same Busemann function}.
\vspace{2mm}

This leads us to describe the Busemann equivalence in a different way.\\
 For quotients $X= G \backslash \tilde{X}$ of Hadamard manifolds, we can characterize Busemann equivalent rays  it in terms  of the dynamics of $G$ on the universal covering  of  $X$. Recall that, by Property (v) after Definition \ref{defiboundary}, the action of $G$ on $\tilde{X}$ extends in a natural way to an action by homeomorphisms on $\partial \tilde{X}$, which is properly discontinuous  on $\mbox{Ord} G$. 

\begin{definition}{$G$-equivalent and weakly  $G$-equivalent  rays.}
\label{defieqG}
{\em 

\noindent Let $\alpha, \beta$ be quasi-rays with origins $a,b$, and lift them to rays  $\widetilde{\alpha}, \widetilde{\beta}$  in $ \tilde{X}$, with origins $\tilde{a},\tilde{b}$. We say that:
\vspace{1mm}

\noindent $\bullet$ $\alpha$ and $\beta$ are  \emph{$G$-equivalent}  ($\alpha \approx_G \beta$) if $\widetilde{\alpha}^+ \in G\widetilde{\beta}^+$; 
\vspace{1mm}

\noindent $\bullet$ {\em $\alpha$ is  weakly $G$-equivalent  to $\beta$} ($\alpha \prec_G \beta$) if there exists
a sequence $g_n \in G$ such that $g_n \widetilde{\beta}^+ \rightarrow \widetilde{\alpha}^+$ and the quasi-rays
 $\alpha_n = \pi [\tilde{a}, g_n \widetilde{\beta}^+]$ have $\Delta(\alpha_n) \rightarrow 0$. \linebreak 
This is equivalent
\footnote{As $\alpha_n= \pi [\tilde{a}, g_n \widetilde{\beta}^+]=\pi [g_n^{-1}\tilde{a}, \widetilde{\beta}^+]$, the excess condition says that $d( g_n^{-1}  \tilde{a} , H^{max}_{\tilde{\beta}^+} (a)) \rightarrow 0$; 
by the formula (\ref{eqbusemann}), this means that
$ B_{\tilde{\beta}} (\tilde{b}, g_n^{-1} \tilde{a}) \rightarrow  \rho \left(  H^{max}_{\tilde{\beta}^+} (b),  H^{max}_{\tilde{\beta}^+} (a)\right) = B_{\beta}(b,a)$.}
 to asking that  there exists a sequence $(g_n)$ such that
\vspace{-3mm}

$$g_n \widetilde{\beta}^+ \rightarrow \widetilde{\alpha}^+ \;\;\; \mbox{ and }\;\;\; B_{\tilde{\beta}} (\tilde{b}, g_n^{-1} \tilde{a}) \rightarrow  B_{\beta}(b,a)$$

\vspace{-1mm}
\noindent where the second condition geometrically means that 
$d( g_n^{-1}  \tilde{a} , H^{max}_{\tilde{\beta}^+} (a)) \rightarrow 0$.

\noindent We say that \emph{ $\alpha$ and $\beta$ are   weakly $G$-equivalent}   ($\alpha \prec_G\succ  \beta$) if   $\alpha \prec_G \beta $ and \nolinebreak  $\beta \prec_G \alpha $.

}
\end{definition}

\noindent Obviously,   {\em $G$-equivalent rays always are weakly asymptotic} (as they admit lifts with common point at infinity); the converse is false in general, as the Example \ref{exladder} in Section \ref{ex} will show.
Further, notice that   {\em $G$-equivalent rays $\alpha, \beta$ define the same Busemann function}; in fact, if $\alpha^+=g\beta^+$,  then according to Theorem \ref{teorbusemann}
\vspace{-5mm}

$$B_{\alpha} (x,y) = \rho (H_{\tilde{\alpha}^+} ^{max}(x),  H_{\tilde{\alpha}^+} ^{max}(y)) 
=\rho (gH_{\tilde{\beta}^+} ^{max}(x),  gH_{\tilde{\beta}^+} ^{max}(y)) 
=B_{\beta} (x,y)$$

\vspace{-2mm}
\noindent but we will see that, in general, {\em two Busemann equivalent rays need not to be $G$-equivalent} (Example \ref{exsymmflute}).
\vspace{2mm}

The interest of the weak $G$-equivalence is explained by the following:

\begin{theorem}
\label{teoreq}
Let $X = G \backslash \widetilde{X}$ be a regular quotient of a Hadamard space. \\
Let $\alpha, \beta$ rays in $X$  with origins $a,b$.  Then:
\vspace{1mm}

\noindent (i) $\alpha \prec \beta$ if and only if $\alpha \prec_G \beta$;
\vspace{1mm}

\noindent (ii) $B_{\alpha} \! = \!  B_{\beta} $ if and only if $\alpha \!  \prec_G\succ  \! \beta$ and
$B_{\alpha} (a,b) = B_{\beta} (a,b)$.
\end{theorem}

\noindent As a corollary, for rays with the same origin $o$, we obtain the Criterium
\ref{critequivalence}.

 \vspace{3mm}
\emph{Proof of Theorem \ref{teoreq}.}  \\
Lift $\alpha, \beta$ to $\widetilde{\alpha}, \widetilde{\beta}$ on  $\widetilde{X}$ with origins $\tilde{a}, \tilde{b}$.
By  Proposition \ref{propcoray}, we know that $\alpha \prec \beta$ if only if $B_{\beta} (a, \alpha(t))=B_{\alpha} (a, \alpha(t) ) =t$ for all $t$.
On the other hand, we have
\vspace{-6mm}

$$ B_{\beta} (a, \alpha(t))   =B_{\beta} (a,b) + B_{\beta} (b, \alpha(t)) =
     B_{\beta} (a,b) + \sup_{g\in G} [ B_{\tilde{\beta}} (\tilde{b}, g \tilde{a}) + B_{\tilde{\beta}} (g \tilde{a}, g \widetilde{\alpha}(t)) ] $$

\vspace{-8mm}

$$ \leq  -B_{\beta} (b,a) + \sup_{g \in G} B_{\tilde{\beta}} (\tilde{b}, g \tilde{a}) + \sup_{g \in G} B_{g^{-1}\widetilde{\beta}} ( \tilde{a},  \widetilde{\alpha}(t))
      \leq d (\tilde{a}, \widetilde{\alpha}(t)) =t $$
 
\vspace{-2mm}
\noindent so $\alpha \prec \beta$ precisely if  there exists a sequence $g_n \in G$ such that
$B_{\tilde{\beta}} (\tilde{b}, g_n \tilde{a}) \rightarrow B_{\beta} (b,a) $ and
$B_{g_n^{-1}\widetilde{\beta}} ( \tilde{a},  \widetilde{\alpha}(t)) \rightarrow   d (\tilde{a}, \widetilde{\alpha}(t)) = B_{\alpha} (a, \alpha(t) )$, that is  $g_n^{-1}\widetilde{\beta}^{+} \rightarrow \widetilde{\alpha}^+$. This shows (i). Part (ii) follows from 
  Theorem \ref{teorcoray}(b).\qed

\subsection{Lower semi-continuity.}
\label{cont}

The behaviour of Busemann functions with respect to the initial directions of quasi-rays is intimately related with the excess. \\
On the one hand, a limit of quasi-minimizing directions does not usually give a direction for which the Busemann function is defined:   
for instance,  if $X=G \backslash \tilde{X}$ and the limit set $LG$ contains at least a Dirichlet point $\zeta$ and a radial point $\xi$, then (as $LG$ is the minimal $G$-invariant closed subset of $\partial \tilde{X}$) there also exists a sequence 
$\zeta_n=g_n \zeta \rightarrow \xi$; the projections $\alpha_n$ on $X$ of rays $[\widetilde{o}, \zeta_n]$ give a family of $G$-equivalent quasi-rays, all defining  the same Busemann function, while the limit curve $\alpha$ is the projection of $[\widetilde{o}, \xi]$, and is a recurrent geodesic for which the Busemann function is not defined.

\noindent Even when the limit curve is a ray or a quasi-ray, with no control of the excess of the family we cannot expect any continuity, as the following example shows:

\begin{example}
\label{exdisc1}
{\em 
\noindent Let $G < Is(\Hy^2)$ be a discrete subgroup generated by two parabolic isometries
$p,q$ with distinct, fixed points $\zeta$, $\xi$, and assume them {\em in Schottky position}, that is: 
$\left( \H^2 - D(<\!p\!>, \widetilde{o}) \right)  \cap \left(   \H^2 - D(<\!q\!>, \widetilde{o}) \right) = \emptyset$, for some $ \widetilde{o} \in \H^2$. \\
For instance, we can take the  group $\G (2)$, generated by  $p(z)= \frac{z}{2z+1}$ and  $q(z)= z+2$ in the Poincar\'e half-plane model, with $\widetilde{o}=i$. In this case, $LG= \partial \Hy^2$  and  $\partial D(G, \widetilde{o})$ consists of two parabolic fixed points $\zeta=0, \xi=\infty$ and  two $G$-equivalent points $\omega=-1$ and $\omega'=1$. The quotient surface $X= G \backslash \Hy^2$ has three cusps corresponding to $\zeta, \xi$ and $\omega'=p(\omega)=q(\omega)$, and only four rays with origin $\widetilde{o}$: the projections $\alpha, \beta, \gamma$ and $\gamma'$ of, respectively, $[\widetilde{o}, \zeta]$, $[\widetilde{o}, \xi]$, 
$[\widetilde{o}, \omega]$ and $[\widetilde{o}, \omega']$, only the last two of which being Busemann-equivalent. \\
By minimality of  $LG$,  there exists a sequence $\zeta_n= g_n \zeta \rightarrow \xi$; then, the projections $\alpha_n$ on $X$ of the rays $[\widetilde{o}, \zeta_n]$ are all $G$-equivalent quasi-rays  (by Corollary \ref{cordictionary}, the $\zeta_n$ being horospherical) which tend to  $\beta$.
However  $B_{\alpha_n}  = B_{\alpha}$ for all $n$, therefore their limit is  $B_{\alpha}$,  while the Busemann function of the limit curve is $B_{\beta} \neq B_{\alpha}$.
}
\end{example}

Notice that in the above examples the excess of the $\alpha_n$ tends to infinity (by Lemma \ref{lemmaexcess}(ii) in the first case, and by direct computation or by the Proposition \ref{propsemicont}  below in Example \ref{exdisc1}).
Keeping control of the length excess yields at least {\em lower semi-continuity} of the Busemann function with respect to the initial directions: 

\begin{proposition}
\label{propsemicont}
{\em Let $X = G \backslash \widetilde{X}$ be a regular quotient of a Hadamard space.
 
\noindent (i) For any sequence of rays  $\alpha_n \rightarrow \alpha$  uniformly on compacts, we have: 
 \vspace{-3mm}
 
\[ \liminf_{n \rightarrow +\infty} B_{\alpha_n} (x,y) \geq  B_{\alpha} (x,y) \]

 \vspace{-1mm} 
\noindent (ii) For any sequence of quasi-rays $\alpha_n \rightarrow \alpha$   with $\Delta(\alpha_n) \rightarrow \Delta (\alpha) + \delta $ we have:
\vspace{-3mm}
 
\[ \liminf_{n \rightarrow +\infty} B_{\alpha_n} (x,y) \geq B_{\alpha} (x,y) - \delta \]

}
\end{proposition}

 \vspace{-2mm} 
{\em Proof.} Part (i) is a particular case of (ii). So let $\widetilde{\alpha}_n, \widetilde{\alpha}$ be lifts of the quasi-rays $\alpha_n, \alpha$ to $\tilde{X}$, with origins $\tilde{a}_n, \tilde{a}$ with $\tilde{a}_n \rightarrow \tilde{a}$, projecting respectively to $a_n, a$. By the cocycle condition, we may assume that $x=a$. By the (\ref{eqbusemanna}) we deduce
\vspace{-3mm}

$$B_{\alpha_n} (a,y) \geq  B_{\tilde{\alpha}_n} (\tilde{a}, g \tilde{y}) - \Delta ( \alpha_n ) - 2d(a, a_n)$$

\vspace{-1mm}
\noindent for all $g \in G$. As $\widetilde{\alpha}_n^+$ tends to $\widetilde{\alpha}^+$ in $\partial \tilde{X}$, we have convergence on compacts of $B_{\tilde{\alpha}_n}$ to $B_{\tilde{\alpha}}$; hence, taking limits for $n \rightarrow \infty$ yields
\vspace{-3mm}

$$ \liminf_{n \rightarrow \infty} B_{\alpha_n} (a,y) \geq  B_{\tilde{\alpha}} (\tilde{a}, g \tilde{y}) - \Delta ( \alpha ) - \delta$$

\vspace{-2mm}
\noindent for all $g$, and we conclude again by using formula (\ref{eqbusemanna}).\qed
\vspace{2mm}

\noindent Lower semi-continuity is the best we can expect, in general,  for the Busemann map:  
in Example  \ref{extwistedflute} we will  produce a case where the strict inequality \linebreak 
 $ B_{\alpha} < \liminf_{n \rightarrow +\infty} B_{\alpha_n}$ holds, for a sequence of rays $\alpha_n \rightarrow \alpha$.

\section{Geometrically finite manifolds}
\label{gf}

 We recall the definition and some  properties of geometrically finite groups.\\
 Let  $G$ be a Kleinian group, that is a discrete, torsionless group of isometries of a negatively curved  simply connected space $\tilde{X}$ with $-a^2 < k(\tilde{X}) \leq -b^2 <0$. \\
Let $\widetilde{C}_G \subset \tilde{X}$ be the convex hull  of  the limit set $LG$;   the quotient  $C_G:=G\backslash \widetilde{C}_G $ is called the {\it Nielsen core} of the manifold $X  =G \backslash \tilde{X}$. The Nielsen core is the relevant subset\footnote{$C_G$ coincides with the smallest closed and convex subset of $X$ containing all the geodesics which meet infinitely often any fixed compact set.}
 of $X$  where the dynamics of geodesics takes place.  \\
    The group  $G$ (equivalently, the manifold $X$) is {\em geometrically finite} if some (any) $\epsilon$-neighbourhood of $C_G$ in $X$  has finite volume.  The simplest examples of geometrically finite manifolds are the {\em lattices}, that is Kleinian groups $G$ such that $vol(G  \backslash \tilde{X} ) < +\infty$. In dimension 2, the class of geometrically finite groups coincides with that of finitely generated Kleinian groups; in dimension $n>2$, geometrically finiteness is a condition strictly stronger than being finitely generated, cp. \cite{apanasov}.  
The following resumes most of the main properties of geometrically finite groups that we will use:  
       
\begin{proposition}[see \cite{bow}]
\label{propgfproperties}
Let $X= G \backslash \tilde{X}$ be a geometrically finite manifold:
\vspace{1mm}

\noindent (a) $LG$ is the union of 
 its  {\it radial} subset  $L^{rad}G$  and of a set $L^{b.par}G= \sqcup_{i=1}^l G\xi_i$ made up of finitely many orbits  of {\em bounded parabolic} fixed points;  this means that each  $\xi \in L^{b.par}G $  is the fixed point of some  maximal parabolic subgroup  ${P}$ of $G$ {\em acting cocompactly on $LG -\xi$}; equivalently, $P$ preserves  every horoball $H_{\xi}$ centered at $\xi$ and acts cocompactly on $\partial H_{\xi} \cap  \widetilde{C}_G$; 
  
 \vspace{1mm}
 \noindent (b) {\em (Margulis' Lemma)} there exist  closed horoballs $H_{\xi_{1}}, \ldots, H_{\xi_{l}}$ centered respectively  at $\xi_{1}, \ldots, \xi_{l}$, such that  $gH_{\xi_{i}} \cap H_{\xi_{j}} = \emptyset$  for  all $1\leq i,j\leq l$ and  all $g \in G - P_i$.

\end{proposition}

\noindent Accordingly, geometrically finite manifolds fall in two classes:
\vspace{1mm}

\noindent $\bullet$  either  $C_G$ is  compact, and   then  $G$ (and $X$) is called  {\em convex-cocompact};
\vspace{1mm}

\noindent $\bullet$  or  $C_G$ is not compact, in which case it can  be decomposed into
 a disjoint union of a compact part $C_0$ and finitely many ``cuspidal  ends'' $C_1, ..., C_l $: each $C_i$ is isometric to the quotient, by  a maximal parabolic group $P_i \subset G$,  of  the intersection between $\widetilde{C}_G \cap H_{\xi_{i}}$, where $H_{\xi_{i}} $ is a horoball preserved by $P_i$ and  centered at ${\xi_{i}}$.

\noindent This  yields a first  topological description of geometrically finite manifolds; for more details on the topology of a horosphere quotient see \cite{belekapo}. In the sequel, we shall always tacitly assume that $X$ is non-compact.
\vspace{2mm}

We will also repeatedly use the following facts:

 \begin{lemma}
\label{lemmagf}
Let $X= G \backslash \tilde{X}$ be a geometrically finite manifold and let $\xi \in LG$ a  bounded parabolic point, fixed by some maximal parabolic subgroup $P< G$:

\noindent (i) $\xi$ is  non-horospherical and universal Dirichlet;

\noindent (ii) there exists a subset $G_{\xi} \subset G$ of representatives of $P\backslash G$ such that  $\xi \not\in \overline{G_{\xi} \tilde{x}}$, for every  $\tilde{x} \in \tilde{X}$.
 \end{lemma}

{\em Proof}.
By Proposition \ref{propgfproperties}(a), we know that $\xi=g \xi_i$ for some $g \in G$, $\xi_i \in P_i$
and that $P=g_iP_ig_i^{-1}$.
Then, consider the family of horoballs $H_{\xi_i}$ given by the Margulis' Lemma, let $H_{\xi}=gH_{\xi_i}$ and  choose a point $\tilde{x}_0 \in \partial  H_{\xi}$, projecting to $x_0 \in X$.
By  Margulis' Lemma, we know that there is no point of the orbit $G\tilde{x}_0$ inside $H_{\xi}$, hence  $H^{max}_{\xi} (x_0)=H_{\xi}$ and $\xi$ is non-horospherical.

\noindent To see (ii), fix a compact fundamental domain $K$ for the action of $G$ on $LG - \xi$: 
then, define the subset $G_{\xi}$ by choosing the identity of $G$ as representative of the class $P$ and, for every $g \in G- P$, a representative $\hat{g}\in Pg$ such that $\hat{g} \xi \in K$.
Since $K$ is compact in $LG-\xi$, it is separated from $\xi$ by an open neighbourhood  $U_{K}$ of $K$ in $\overline{X}$, $\xi \not\in U_{K}$.
Now, as $\xi$ is universal Dirichlet, for every fixed $\tilde{x}$ we can find $g \in G$ such that
$\xi \in \partial D ( G, g \tilde{x})$; by construction, the orbit $G_{\xi} \xi$ accumulates to $K$ and, as the Dirichlet domain is locally finite, the domain  $D ( G, g \tilde{x})$ too.
Since $d( \tilde{x}, g \tilde{x}) < \infty$, we also deduce that the subset $G_{\xi}  \tilde{x}$  is included  (up to a finite subset) in $U_K$; this shows that $\xi \not\in \overline{G_{\xi} \tilde{x}}$.\qed

\begin{proposition}
\label{propprerays}
Let $X= G \backslash \tilde{X}$ be a geometrically finite manifold: then, 
every quasi-ray of $X$ is a pre-ray.
\end{proposition}

\emph{Proof}. 
Let $\alpha$ be a quasi-ray of $X$  with origin $a$, and lift it to a ray $\widetilde{\alpha}$ of  $\tilde{X}$, with origin $\tilde{a}$.
Assume that $\alpha$ is not a pre-ray: then, by Property \ref{lemmaexcess}(i)  we would have a positive, strictly decreasing sequence $\Delta_n=\Delta( \alpha |_{[t_n, +\infty)} )$, tending to zero, for some  $t_n \rightarrow +\infty$.
Since $\widetilde{\alpha}^+$ is a non-horospherical point, it is either ordinary or bounded parabolic; anyway, it is a universal Dirichlet point by Lemma \ref{lemmagf}, so for each $n$ we can find $g_n$   such that $g_n^{-1}  \widetilde{\alpha} (t_n) \in \partial H^{max}_{\tilde{\alpha}^+} (\alpha (t_n))$. \linebreak
Let $P$ be the maximal parabolic subgroup fixing $\xi=\tilde{\alpha}^+$, and  let $\hat{g}_n=p_ng_n$ be the representative of $g_n \in G_{\xi}$ given by Lemma \ref{lemmagf}, for $p_n \in P$. 
We have:
\vspace{-5mm}

$$\Delta (\alpha) \geq B_{\tilde{\alpha}} ( \widetilde{a}, g_n^{-1} \tilde{a}) =
  B_{\tilde{\alpha}} ( \widetilde{a}, \hat{g}_n^{-1} \tilde{a}) =
B_{\tilde{\alpha}} ( \tilde{a},  \widetilde{\alpha} (t_n)) +
B_{\tilde{\alpha}} (  \widetilde{\alpha} (t_n), \hat{g}_n^{-1} \widetilde{\alpha} (t_n)) +$$

\vspace{-6mm}

$$ + B_{\tilde{\alpha}} (  \hat{g}_n^{-1} \widetilde{\alpha} (t_n), \hat{g}_n^{-1}  \tilde{a})\geq  t_n + \Delta_n -  B_{\hat{g}_n\tilde{\alpha}} (\tilde{a},   \widetilde{\alpha} (t_n))$$

\vspace{-1mm}
\noindent which, since the excess of $\alpha$ is finite,  shows that $\hat{g}_n\tilde{\alpha}^+ \rightarrow \tilde{\alpha}^+$ necessarily, for $n\rightarrow \infty$. 
By the locally finiteness of the Dirichlet domain, we deduce that  $\hat{g}_n\tilde{a}  \rightarrow \tilde{\alpha}^+$ too, which contradicts (ii) of Lemma \ref{lemmagf}.\qed

  \pagebreak 
 For geometrically finite manifolds, the equivalence problem is answered by:
  
\begin{proposition}
\label{propeq}
Let $X= G \backslash \tilde{X}$ be a geometrically finite manifold, and let $\alpha, \beta$ rays. The following conditions are equivalent:
\vspace{1mm}

\noindent
(a) $\; B_{\alpha}=B_{\beta}$ \hspace{10mm}
 (b) $\;\alpha \approx_G \beta$ \hspace{10mm}
(c) $\;\alpha \prec \beta$\hspace{10mm}
 (d) $\; d_{\infty}(\alpha, \beta) < \infty$
\end{proposition}

\emph{Proof}. 
Let $a,b$ the origins of the two rays $\alpha, \beta$, and let $\widetilde{\alpha}$ and $\widetilde{\beta}$ the lifts of $\alpha, \beta$ to $\tilde{X}$, with origins $\tilde{a}, \tilde{b}$ respectively. 
Now assume that $\alpha \prec \beta$. 
Consider the quasi-ray $\beta'$ which is the projection of $\tilde{\beta}'=[\tilde{a}, \tilde{\beta}^+]$ to $X$,  and fix a  $t_0>0$.\linebreak
Since $\alpha \prec \beta \approx_G \beta'$ we have, by Proposition \ref{propcoray} and Theorem \ref{teoreq}
\vspace{-3mm}

$$B_{\beta'}(a,\alpha(t_0)) = B_{\beta} (a,\alpha(t_0)) = d(a, \alpha(t_0))  = t_0.$$

\vspace{-1mm}
\noindent As $G$ is geometrically finite, $\widetilde{\beta}^+$ is  universal Dirichlet and there exists $g_0$ such that $g_0\widetilde{\alpha}(t_0) \in \partial H^{max}_{\tilde{\beta}^+}(\alpha(t_0))$.
Then,  by Theorem \ref{teorbusemann} and the Excess Lemma \ref{lemmaexcess}
\vspace{-3mm}

$$t_0=B_{\beta'}(a,\alpha(t_0))
=  B_{\tilde{\beta}'}(\tilde{a},g_0\tilde{a} ) +  B_{\tilde{\beta}'}(g_0\tilde{a},g_0\widetilde{\alpha}(t_0)) - \Delta(\beta') \leq$$

\vspace{-5mm}

$$\leq  B_{g_0^{-1}\tilde{\beta}'}(\tilde{a},\widetilde{\alpha}(t_0))  \leq d(\tilde{a}, \widetilde{\alpha}(t_0))=t_0$$

\vspace{-1mm}
\noindent Then $B_{g_0^{-1}\tilde{\beta}'}(\tilde{a},\widetilde{\alpha}(t_0) ) =d(\tilde{a}, \widetilde{\alpha}(t_0))$, hence 
$g_0^{-1}\widetilde{\beta}^+ =g_0^{-1}\widetilde{\beta}'^+ = \widetilde{\alpha}^+$. Therefore $\alpha \approx_G \beta$, which implies $B_{\alpha} = B_{\beta}$. As (a) implies (c), this shows that the first three conditions are equivalent. 
To conclude, let us show that (d) and  (b) are equivalent. We already remarked that $G$-equivalence implies  asymptoticity. So, assume now that $d_{\infty} (\alpha, \beta) < +\infty$. 
Up to replacing $\beta$ with the $G$-equivalent quasi-ray $\beta'$  defined above, which still has $d_{\infty} (\alpha, \beta') \leq M< +\infty$, we can assume that their lifts $\widetilde{\alpha}$ and $\widetilde{\beta}$ have the same origin $\tilde{a}$.
Then, let $t_k, t'_k \rightarrow +\infty$ and $g_k \in G$ such that $d(\alpha(t_k), \beta(t'_k)) = d(\widetilde{\alpha}(t_k), g_k \widetilde{\beta}(t'_k)) \leq M$; this implies that $g_k \beta^+ \rightarrow \alpha^+$. Now, if the $g_k$'s form a finite set, then $g_k \widetilde{\beta}^+ = \widetilde{\alpha}^+$ for some $k$, and the rays are $G$-equivalent.
Otherwise, since $G$ acts discontinuously on $\partial \tilde{X} - LG$, we deduce that $\widetilde{\alpha}^+ \in LG$; moreover, as $\widetilde{\alpha}^+$ is a Dirichlet point, it necessarily is  a bounded parabolic point of $G$. We deduce analogously that $\widetilde{\beta}^+$ is parabolic. But now, if $\tilde{\beta}^+ \not\in G\tilde{\alpha}^+$, Margulis' Lemma yields horoballs $H_{\tilde{\alpha}^+},  H_{\tilde{\beta}^+}$,  respectively containing
$\widetilde{\alpha}(t_k)$ and $\widetilde{\beta}(t_k)$ for $k\gg 0$, such that $H_{\tilde{\alpha}^+} \cap  gH_{\tilde{\beta}^+} = \emptyset$  for all $g \in G$.
Then,  $d( \widetilde{\alpha}(t_k), g_k\widetilde{\beta}(t_k) ) \geq
d( \widetilde{\alpha} (t_k), H_{\tilde{\alpha}^+}) \rightarrow  + \infty$,
which contradicts our assumption.\qed

\begin{proposition} 
\label{propsur}
Let $X= G \backslash \tilde{X}$ be a geometrically finite manifold. \\ For any $o \! \in\! X$ the Busemann map $B_o \!: {\cal R}_o X \!\rightarrow \! \partial X$ is surjective, i.e. \nolinebreak ${\cal B}X \!=\! \partial X$. \\
Namely, let $(x_n)$  be a sequence of points converging to a horofunction  $B_{(x_n)}$. \linebreak If $ \tilde{x}_n$ are lifts of the $x_n$ in a Dirichlet domain $D(G,\tilde{o})$, accumulating to some $\xi \in \partial D(G,o)$, then $B_{(x_n)} = B_{\alpha}$ where $\alpha$ is the ray projection of  $[\tilde{o},\xi]$ to $X$.
\end{proposition}

{\em Proof.}  First notice that we have
$d(o, x_n) - d(x_n,x) \geq d(\tilde{o}, \tilde{x}_n) -d(\tilde{x}_n, g\tilde{x})$ for every $g$  and, by taking limits,  we get $B_{(x_n)} (o,x) \geq B_{\xi} (\tilde{o},g\tilde{x})$,  as the $\tilde{x}_n$ accumulate to $\xi$;  therefore  $B_{(x_n)} (o,x) \geq \sup_g B_{\xi} (\tilde{o},g\tilde{x})= B_{\alpha} (o,x)$, by (\ref{eqbusemannb}).  \\
To show the converse inequality, let $x$ be fixed and for each $n$ choose $g_n$  such that $d(x, x_n) = d (g_n \tilde{x}, \tilde{x}_n)$. We will show that there exists  $\hat{g} \in G$ such that 
\vspace{-2mm}

\begin{equation}
\label{limitg}
  d (g_n \tilde{x}, \tilde{x}_n) - d( \hat{g}  \tilde{x}, \tilde{x}_n) \longrightarrow 0 \;\;\;\mbox{ as } n \rightarrow \infty
\end{equation} 

\noindent up to a subsequence; then, from this we will deduce that
\vspace{-3mm}

$$[ d(o,x_n) - d(x_n, x) ] -  [ d(\tilde{o},  \tilde{x}_n) - d( \tilde{x}_n, \hat{g} \tilde{x})  ] \; \longrightarrow 0\;$$

\vspace{-1mm}
\noindent and, as the first summand tends to $B_{(x_n)} (o,x)$ and the second to  $B_{\xi} (\tilde{o}, \hat{g}\tilde{x})$, we can conclude that
$ B_{(x_n)} (o,x) = B_{\xi} (\tilde{o}, \hat{g}\tilde{x}) \leq  \sup_g B_{\xi} (\tilde{o},g\tilde{x})= B_{\alpha} (o,x)$.

\noindent Let us then show (\ref{limitg}). Notice that this is evident  when the set of the $g_n$ is finite.
So, assume  that the set is infinite;  then $g_n \tilde{x}$ accumulates to some limit point $\eta$. \linebreak
If $\eta \neq \xi$, let $\vartheta_0= \widehat{ \xi \; \tilde \!\! o \eta} >0$; then, by comparison geometry,  there exists  $c=c(\vartheta_0)$ (also depending on the upper bounds of the sectional curvature of $\tilde{X}$) such that  for $n \gg0$ 
\vspace{-6mm}

$$d( g_n \tilde{x}, \tilde{x}_n) \sim_{c(\vartheta_0)}  d(g_n \tilde{x}, \tilde{o}) + d(\tilde{o},\tilde{x}_n)  $$

\vspace{-1mm}
\noindent but, as  $d(g_n \tilde{x}, \tilde{o}) \rightarrow + \infty$,  this contradics the fact that $d(\tilde{o}, \tilde{x}_n) - d(\tilde{x}_n, g_n \tilde{x})$ converges. Therefore $g_n \tilde{x} \rightarrow \xi \in LG \cap \partial D(G,\tilde{o})$, and $\xi$ necessarily is a bounded parabolic point.
Then, let $P$ be the maximal parabolic subgroup fixing $\xi$, and  let $\hat{g}_n=p_ng_n$ be the representative of $g_n$ in the subset $G_{\xi}$ given by Lemma \ref{lemmagf}, for $p_n \in P$.
We again have that $p_n x_n \rightarrow \xi$ up to a subsequence; in fact,   $x_n$ tends to $\xi$ within $D(G, \tilde{o})$, so either  the $p_n$'s form a finite set and $p_n x_n = p x_n \rightarrow p \xi=\xi$, or  the whole  $p_n D(G, \tilde{o})$ converges to $\xi$  (the Dirichlet domain being locally finite).
We now infer that the set of $\hat{g}_n$ is finite: otherwise, the points $ \hat{g}_n \tilde{x} = p_n g_n  \tilde{x}$ would accumulate to some $\eta$ different from $\xi$ (by Lemma \ref{lemmagf}); and the same comparison argument as above would give
\vspace{-6mm}

$$d(x, x_n) =  d(g_n \tilde{x}, \tilde{x}_n) = d (\hat{g}_n \tilde{x}, p_n \tilde{x}_n  )
 \sim_{c(\vartheta_0)}  d(\hat{g}_n \tilde{x}, \tilde{x}) + d(\tilde{x},p_n\tilde{x}_n)  
 \gg   d(x, x_n) $$

\vspace{-1mm}
\noindent for $n$ large enough, which is a contradiction. Thus, the set of $\hat{g}_n$ is finite,  and we may assume that $g_n = \hat{g}$ definitely. Now
\vspace{-3mm}

\begin{equation}
\label{eqsomma}
 [ d(\tilde{o}, \tilde{x}_n) - d(\tilde{o}, p_n \tilde{x}_n) ] + 
     [ d( \tilde{x}_n, p_n^{-1}\hat{g} \tilde{x}) - d( \tilde{x}_n, \hat{g} \tilde{x}) ] 
\end{equation}

\vspace{-5mm}

\begin{equation}
\label{eqsomma2}
\hspace{1cm} =  [ d(\tilde{o}, \tilde{x}_n) - d( \tilde{x}_n, \hat{g} \tilde{x}) ] -
    [  d(\tilde{o}, p_n \tilde{x}_n) - d( p_n \tilde{x}_n, \hat{g} \tilde{x})  ] \rightarrow  0
\end{equation}
as we know that both $\tilde{x}_n$ and $p_n \tilde{x}_n$ tend to $\xi$, so both terms in   (\ref{eqsomma2}) tend to
$B_{\xi} (\tilde{o}, \hat{g} \tilde{x})$.
The first summand  $ [ d(\tilde{o}, \tilde{x}_n) - d(\tilde{o}, p_n \tilde{x}_n) ] $ in (\ref{eqsomma}) is nonpositive since the $\tilde{x}_n$ belong to $D(G, \tilde{o})$;  the second summand in (\ref{eqsomma}) also is nonpositive, as
\vspace{-3mm}

$$d( \tilde{x}_n,  p_n^{-1}\hat{g} \tilde{x}) = d( \tilde{x}_n, g_n \tilde{x}) \leq  d( \tilde{x}_n, g \tilde{x})  \;\;\;\; \forall g \in G$$

\vspace{-1mm}
\noindent by assumption; therefore by (\ref{eqsomma2}) we deduce that $  d( \tilde{x}_n,  g_n \tilde{x}) - d( \tilde{x}_n, \hat{g} \tilde{x} ) \rightarrow 0$
which proves (\ref{limitg}) and concludes the proof.\qed \\

\vspace{3mm}
For the next result, we need to recall the  Gromov-Bourdon metric on   $\partial \tilde{X}$.
This is a family of metrics  indexed by the choice of a base point $\tilde{o} \in \tilde{X}$:
\vspace{-3mm}

$$ D_{\tilde{o}} (\eta, \xi) = e^{- \frac12 | B_{\eta}  (\tilde{o}, \tilde{x}) +  B_{\xi}  (\tilde{o}, \tilde{x}) |} 
\;\;\; \mbox{ for any } \tilde{x} \in [\eta, \xi]$$

\vspace{-1mm}
\noindent
The exponent  corresponds to minus the length of the finite geodesic segment  cut on    $[\eta, \xi]$ by the horospheres $H_{\eta} (\tilde{o})$, $H_{\xi} (\tilde{o})$.
The fundamental property of these metrics is that any isometry of $ \tilde{X}$ acts by conformal homeomorphisms on  $\partial  \tilde{X}$ with respect to them; moreover, the conformal coefficient can be easily expressed in terms of the Busemann function \cite{bourdon}:
\vspace{-5mm}

\begin{equation}
\label{bourdon}
D_{\tilde{o}} (g\eta, g\xi) = \sqrt{g' (\eta)} \sqrt{g' (\xi)}   D_{\tilde{o}} (\eta, \xi) \;\;\; \mbox{ where } g' (\zeta) =  e^{ B_{\zeta}  (\tilde{o}, g^{-1} \tilde{o})}
\end{equation}

\begin{proposition}
\label{propcont}
Let $X= G \backslash \tilde{X}$ be a geometrically finite manifold, and let $\alpha_n$ be a sequence of rays converging to $\alpha$. Then, $B_{\alpha_n} (x,y) \rightarrow B_{\alpha} (x,y)$ uniformly on compacts.
\end{proposition}

{\em Proof.} Notice that the limit curve $\alpha$ still is a ray by Lemma \ref{lemmaexcess}. Also, notice that, if $a$ is the origin of $\alpha$, by the cocycle condition it is enough to show that $B_{\alpha_n} (a,x)$ converges uniformly on compacts to $B_{\alpha} (a,x)$.
Then, let $\widetilde{\alpha}$, $\widetilde{\alpha}_n$ be rays of $\tilde{X}$
with origins $\tilde{a}, \tilde{a}_n$ 
projecting  respectively to $\alpha$ and $\alpha_n$, and let $\epsilon_n= d(\tilde{a}, \tilde{a}_n)\rightarrow 0$.
Now choose any  point $x \in X$. Since $G$ is geometrically finite, $\alpha^{+}$ and $\alpha_n^{+}$ are either ordinary or bounded parabolic points; anyway, they are universal Dirichlet, so let $\tilde{x}$ and $g_n \tilde{x}$ be  lifts of $x$ such that
\vspace{-3mm}

$$B_{\alpha} (a,x) = B_{\tilde{\alpha}} (\tilde{a},\tilde{x}), \;\;\;\; 
B_{\alpha_n} (a_n,x) =B_{\tilde{\alpha}_n} (\tilde{a}_n,g_n\tilde{x})$$

\vspace{-1mm}
\noindent by Theorem \ref{teorbusemann}.
If $\widetilde{\alpha}^{+}$ is parabolic, let $P$ be its maximal parabolic subgroup and let $\hat{g}_n=p_ng_n$ be the representative of $g_n$ modulo $P$ given by Lemma \ref{lemmagf}, with $p_n \in P$; if $\widetilde{\alpha}^{+}$ is ordinary, just set $\hat{g}_n=g_n$ and $p_n=id$.
Then, consider the set $F$ of all the $\hat{g}_n$'s: we claim that $F$ is finite.
In fact, first notice that
\vspace{-3mm}

$$D_{\tilde{a}}(p_n \widetilde{\alpha}_n^+, \widetilde{\alpha}^{+}) = \sqrt{p_n' (\widetilde{\alpha}_n^+)} D_{\tilde{a}}( \widetilde{\alpha}_n^+, \widetilde{\alpha}^{+}) 
\leq e^{2\epsilon_n} D_{\tilde{a}}( \widetilde{\alpha}_n^+, \widetilde{\alpha}^{+})$$

\vspace{-1mm}
\noindent as $B_{\tilde{\alpha}_n^+} (\tilde{a}, p_n^{-1}\tilde{a})\leq 2 \epsilon_n$,
$\alpha_n$ being a ray from $a_n$ with $d(a,a_n)=\epsilon_n$; therefore, we deduce that $p_n \widetilde{\alpha}_n^+ \rightarrow \widetilde{\alpha}^{+}$.
Moreover, we have
\vspace{-5mm}

\begin{equation}
\label{eqfinite}
-d(a, x)  \!\leq \! B_{\alpha_n} (a,x)  \!= \! B_{\tilde{\alpha}_n} (\tilde{a},p_n^{-1}\tilde{a}) \!+\!  B_{\tilde{\alpha}_n} (p_n^{-1}\tilde{a},g_n\tilde{x}) \apprle_{2\epsilon_n} \!\! B_{p_n\tilde{\alpha}_n} (\tilde{a},\hat{g}_n\tilde{x})
\end{equation}

\vspace{-1mm}
\noindent If $F$ is infinite,  we deduce $\hat{g}_n\tilde{x} \rightarrow \xi \neq \widetilde{\alpha}^{+} $ by Lemma \ref{lemmagf}, so 
$B_{p_n\tilde{\alpha}_n} (\tilde{a},\hat{g}_n\tilde{x}) \rightarrow -\infty$, contradicting (\ref{eqfinite}).
So, $F$ is finite and we may assume that $\hat{g}_n = \hat{g}$ definitely.
But then, passing to limits in (\ref{eqfinite}) we get
\vspace{-3mm}

$$\lim_{n \rightarrow + \infty} B_{\alpha_n} (a,x) \leq   B_{\tilde{\alpha}} (\tilde{a},\hat{g}\tilde{x})
\leq B_{\alpha} (a,x)\; . $$

\vspace{-1mm}
\noindent By the lower semi-continuity (Proposition \ref{propsemicont}) we deduce that $B_{\alpha_n} (a,x)$ converge pointwise to $B_{\alpha} (a,x)$; but as $B_{\alpha_n} (a,x) $ are a family of $1$-Lipschitz functions of $x$, this implies uniform convergence on compacts.\qed

\begin{corollary}
\label{corgf}
Let $X= G \backslash \tilde{X}$ be a geometrically finite, $n$-dimensional manifold. \\ 
For any  $\tilde{o} \in \tilde{X}$ projecting to $o \in X$, the horoboundary $\partial X$ of $X$ is homeomorphic to
\vspace{-5mm}

\begin{equation}
\label{eqeq}
{\cal R}_o(X) /_{(Busemann \; eq.)} \cong G \backslash \partial D(G, \tilde{o})
\end{equation}

\noindent and the horofunction compactification of $X$ is $\overline{X} \cong G \backslash \overline{D(G, \tilde{o})}$. \\
If $n=2$ or $G$ has no parabolic subgroups, then $\overline{X}$ is a topological manifold with boundary. 
If $n \geq 3$ and $G$ has parabolic subgroups, then $\overline{X}$ is a topological manifold with boundary with a finite  number of conical singularities, each corresponding to a conjugate class of maximal parabolic subgroups of $G$.
\end{corollary}

\noindent Here, we call  {\em conical singularity} a point $\xi$ with a neighbourhood homeomorphic to the cone over some topological manifold  (with or without boundary) \nolinebreak $Y$:
\vspace{-2mm}

$$C(Y, \xi)= \left( Y \times [0,1] \right) / _{(y,1)=\xi,\; y \in Y}$$

\vspace{-1mm}
\noindent and we say that $\overline{X}$ is a {\em topological manifold with conical singularities} if $\overline{X}$ has a discrete subset  $S=\{\xi_k\}$ of conical singularities such that $\overline{X} - S$ is a usual topological manifold (with or without boundary).

\vspace{3mm}
{\em Proof.} By the Property \ref{lemmaexcess}(ii), for any $o\in X$ the set of rays from $o$ can be topologically identified to a {\em closed} subset of the tangent sphere $S_oX$ at $o$, hence it is compact. Then, by Propositions \ref{propsur} \& \ref{propcont} we deduce that  the restriction of the Busemann map $B_o: {\cal R}_o (X) /_{(Busemann \; eq.)} \rightarrow \partial X$
  is a homeomorphism.
Moreover, the set of rays of $X$ with origin $o$ consists of all projections of half-geodesics from $\tilde{o}$ in $\tilde{X}$ staying in the Dirichlet domain, i.e. whose boundary points belong to $\partial D(G,\tilde{o})$. Since by Proposition \ref{propeq} the Busemann equivalence is the same as $G$-equivalence, this establish the bijection (\ref{eqeq}). Notice that this is a homeomorphism as the uniform topology on ${\cal R}_o (X)$ corresponds to the sphere topology on (the subset of minimizing directions of) $S_o X$. Then, as $ G \backslash D(G, \tilde{o}) \cong X$,  the map $b$ of Section \ref{busemannfunctions} establishes the homeomorphism $G \backslash \overline{D(G, \tilde{o})} \cong \overline{X}$. 
Let us now precise the structure of $\overline{X}$ at its boundary points. \\
We know that $\partial D(G,\tilde{o})$ is made up of ordinary points of   $\mbox{Ord} G$ and finitely many orbits of bounded parabolic points $\xi_k$; let  $\partial_{ord} D(G,\tilde{o})$ the subset of ordinary points on the trace of the Dirichlet domain. Every ordinary point $\xi \in \mbox{Ord} G$ has a neighbourhood homeomorphic to a neighbourhood of a boundary point of the closed, unitary Euclidean ball in $T_oX$ centered at $0$, and  the action of $G$ on $\mbox{Ord} G$ is proper.
So, the space
\vspace{-3mm}

$$X'=  G \backslash ( \tilde{X} \cup \mbox{Ord} G ) =   G \backslash \left[ D(G,\tilde{o}) \cup \partial_{ord} D(G,\tilde{o})  \right]$$

\vspace{-1mm}
\noindent  has a structure of ordinary topological manifold with boundary. This structure coincides with the uniform topology of the horofunction compactification, as a sequence $(x_n)$ in $ D(G,\tilde{o})$ tend to an ordinary point $\xi$ if and only if $b_{x_n} \rightarrow B_{\xi}$, by Proposition \ref{propsur}.
Now, $X'$ has a finite number of ends $E_k$, corresponding to the classes modulo $G$ of the bounded parabolic points $\xi_k$; we will use the description of  such ends due to Bowditch, to figure out their horofunction compactification.
Let $P_k$ be the maximal parabolic subgroup associated with $\xi_k$, let  $H_{k}$ some  horosphere centered at $\xi_k$, with quotient $Y_k=P_k \backslash H_{\xi_k}$, and let  $X_k=P_k \backslash \tilde{X}$. $X_k$ is a geometrically finite manifold, with one orbit of parabolic points corresponding to $\xi_k$, and the manifold with boundary 
\vspace{-3mm}

$$X'_k= P_k \backslash ( \tilde{X} \cup \mbox{Ord} P_k ) = 
 P_k \backslash \left[ \overline{D(P_k,\tilde{o})} - \xi_k  \right]$$

\vspace{-1mm}
\noindent  has one end {\em isometric to the end $E_k$}, cp. \cite{bow}; topologically, $X'_k=Y_k \times [0, \infty)$.
By \cite{belekapo}, $Y_k$ is a vector bundle over a compact manifold $M_k$, so let ${\cal D}(Y_k)$ and ${\cal S}(Y_k)$ the associated closed disk and sphere bundles. 
The horofunction compactification of the end $E_k$, by (\ref{eqeq}), has just one point  at infinity corresponding to $\xi_k$, 
and  is homeomorphic to
\vspace{-7mm}

 $$C({\cal T}(Y_k), \xi_k) = \frac{{\cal D}(Y_k) \times [0,\infty]}{ \left( {\cal S}(Y_k) \! \times \! [0,\infty] \right) \cup \left( {\cal D}(Y_k)  \! \times \! \{ \infty \} \right)} 
 	= \frac{ {\cal T}(Y_k) \times [0,\infty]   }{ {\cal T}(Y_k)  \times \{ \infty \} },$$

\vspace{-1mm}
\noindent the cone over the Thom space ${\cal T}(Y_k)= {\cal D}(Y_k) / {\cal S}(Y_k)$ of   $Y_k$, with vertex $\xi_k$; 
actually, every sequence of points diverging in the end yields the same horofunction (the Busemann function of the projection to $X_k$ of  $[\tilde{o}, \xi_k]$, by Proposition \nolinebreak \ref{propsur}).  
Notice that, on each fiber of $Y_k$ over $m \in M_k$,  the space
$\frac{{\cal D}_m (Y_k) \times [0,\infty]}{ \left( {\cal S}_m(Y_k) \! \times \! [0,\infty] \right) \cup \left( {\cal D}_m(Y_k)  \! \times \! \{ \infty \} \right)}$
is homeomorphic to the cone  $C({\cal D}_m (Y_k), \xi_{m}(k) )$ with  base ${\cal D}_m (Y_k)$ and  vertex
 $\xi_{m}(k)$;
it follows that 
 \vspace{-3mm}
 
$$ C({\cal T}(Y_k), \xi_k) \cong  \frac{ \bigcup_{m \in M_k} C({\cal D}_m (Y_k), \xi_{m}(k)) }{ \cup_{m \in M_k} \; \xi_{m}(k) }
\cong \frac{ {\cal D}(Y_k) \times [0,\infty]   }{ {\cal D}(Y_k)  \times \{ \infty \} } = C({\cal D}(Y_k), \xi_k) $$

\vspace{-1mm}
\noindent is homeomorphic to the cone over the closed manifold (with boundary) ${\cal D}(Y_k)$.
Clearly ${\cal T}(Y_k)= {\cal D}(Y_k)=Y_k =M_k$  if $n=2$, and in this case $C({\cal D}(Y_k),\xi_k)$ is a closed topological disk;
on the other hand, in dimension $n\geq 3$ this cone is always singular at $\xi_k$ (since $Y_k$ is not simply connected, the subset $C({\cal D}(Y_k), \xi_k)- \xi_k$ is not locally simply connected).\qed

\begin{examples}  The horofunction compactification of an unbounded cusp
\label{extomate}
{\em 

\noindent (i) Let $X =P \backslash \H^3$ where $P$ is generated by a parabolic isometry $p$ with fixed point $\xi$. In the Poincar\'e half-space model, assume that $\xi$ is the point at infinity, fix some horosphere $H_{\xi}$ and choose a origin  $\tilde{o}$. The Dirichlet domain $D(P,\tilde{o})$ is an infinite vertical corridor, with  parallel vertical walls $W_1,W_2$ paired by $p$. \linebreak
 $X$ is homeomorphic to an open cylindrical shell, which is the product of the horosphere quotient $ Y=P \backslash H_{\tilde{\xi}}=\mbox{Cyl} $ (a flat infinite cylinder) with \nolinebreak $\mathbb{R}^{\ast}_+:$
\vspace{-3mm}

$$X=P \backslash D(P,\tilde{o})
\cong  \mbox{Cyl} \times (0, \infty).$$

\vspace{-1mm}
\noindent We may take $\mbox{Cyl} \cong S^1 \times (-1,1)$ with closure $\overline{\mbox{Cyl}} = S^1 \times [-1,1]$ .
Then, the manifold $X'$ is
\vspace{-2mm}

$$X'= P \backslash \left[ \Hy^3 \cup \mbox{Ord} P \right] = P \backslash \left[ \overline{D(P,\tilde{o})} - \xi \right] \cong \mbox{Cyl} \times [0, \infty)$$

\vspace{-1mm}
\noindent the end of which corresponds to a neighbourhood of the bases $B^{\pm} = S^1 \times \{ \pm 1\}$ of the cylinder and of the internal boundary  $\mbox{Cyl}^{\infty}=\mbox{Cyl} \times \{ \infty \}$ of the shell \linebreak (a solid hourglass).
The horofunction compactification is 
\vspace{-3mm}

 $$\overline{X} \cong \overline{\mbox{Cyl}} \times [0, \infty]  \; /_{ (B^+=B^-=Cyl^{^{\infty}})}$$

\vspace{-1mm}
\noindent that is, a spindle solid torus, whose center corresponds to the unique singular point at infinity of the compactification.
\vspace{2mm}

\noindent (ii) Let $X =G \backslash \H^3$ where  $G=<p,h>$ is the free group generated by a parabolic isometry $p$ and a hyperbolic isometry $h$ in Schottky position. In this case, the Dirichlet domain  is the same vertical corridor as above, minus two hemispherical caps (the attractive and repulsive domains of $h$), and the horofunction compactification is  the above spindle solid torus    with a solid handle attached.
}
\end{examples}

\section{Examples}
\label{ex}

\begin{figure}[h]
\hspace{5mm}
\includegraphics[width=110mm, height=27mm]{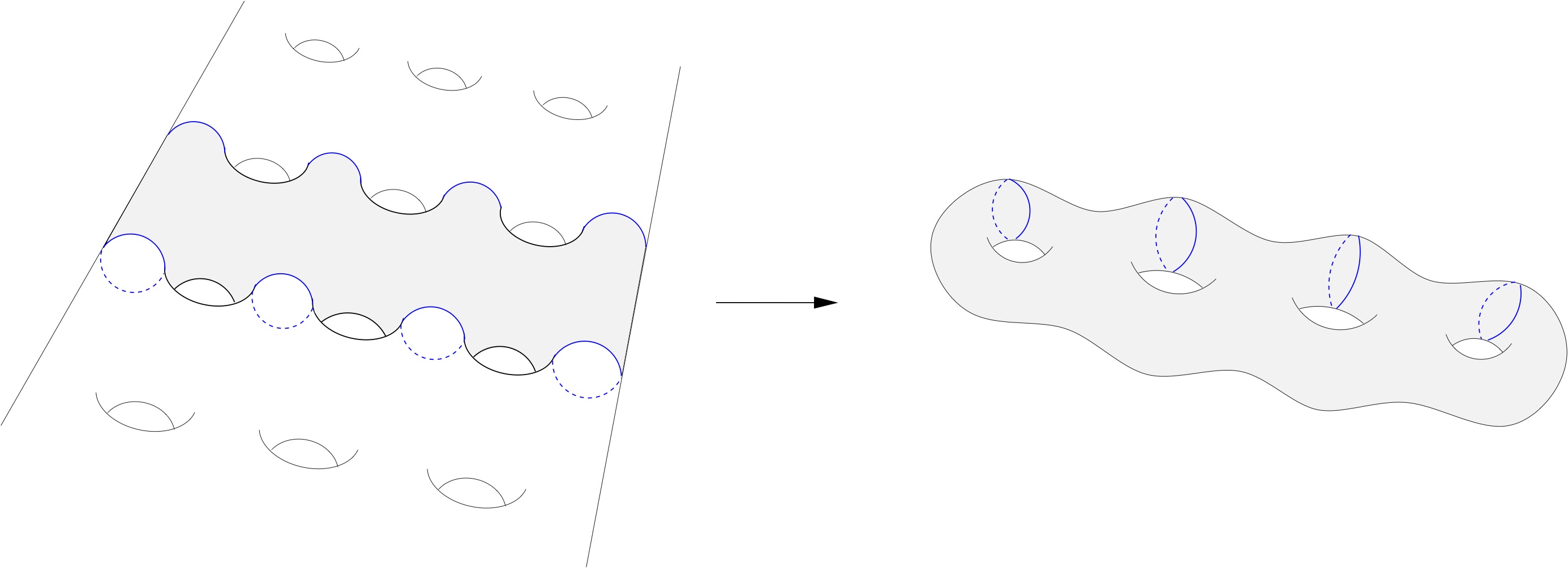}
\caption{Construction of ladders}
\label{figconstructionladder}
\end{figure}

We present in this section some examples of two basic classes of complete, non-geometrically finite hyperbolic surfaces presenting the pathologies described in the introduction  (Theorems \ref{teordifferentrelations}, \ref{teornotclosed}, \ref{teornonsur}, \ref{teornoncont}):

\vspace{1mm}
\noindent 
$\bullet$ {\sc Hyperbolic Ladders}: these are $\mathbb{Z}$-coverings of a hyperbolic closed surface $\Sigma_g$ of genus $g \geq 2$,  obtained by infinitely many copies of the base surface $\Sigma_g$ cut along  $g$ simple, non-intersecting closed geodesics of a fundamental system, glued along the corresponding boundaries, cp. Figure \ref{figconstructionladder};

\noindent 
$\bullet$ {\sc Hyperbolic Flutes}: these are, topologically, spheres with infinitely many punctures $e_i$ accumulating to one limit puncture $e$; the surface thus has one end for each puncture $e_i$ (called its {\em finite ends}), and an end corresponding to $e$, the {\em infinite end} of the flute. Geometrically, each  end $e_i$ different from $e$ must be either a {\em cusp} (the quotient of a horoball $H_{\xi}$ of $\Hy^2$ by a parabolic subgroup $P_{\xi}$ fixing the center $\xi$ of $H_{\xi}$) or a {\em funnel} (the  quotient of a half-plane of $\Hy^2$ by an infinite cyclic group of hyperbolic isometries).
\vspace{2mm}

 We obtain workable models of flutes via infinitely generated {\em Schottky groups}.
Define the  {\em attractive} and {\em repulsive} domains $A(g, \tilde{o})$,   $A(g^{-1}, \tilde{o})$ of a parabolic or hyperbolic isometry $g$, with respect to some point $\tilde{o} \in \Hy^2$, respectively as
\vspace{-3mm}
 
$$A(g^{\pm 1}, \tilde{o}) = \{ x \in \Hy^2 \; | \;   d(x,\tilde{o}) \geq  d(x,g^{\pm 1}\tilde{o})   \}$$

\vspace{-1mm}
\noindent We say that $G$ is an infinitely generated Schottky group if it is generated by countable many  hyperbolic isometries $S=(g_n)$,  in {\em Schottky position with respect to some $\tilde{o} \in \Hy^2$}, that is:
 $A(g_n^{\epsilon}, \tilde{o}) \cap A(g_m^{\epsilon'}, \tilde{o}) = \emptyset
\;\;\; \forall n, m \mbox{ and } \forall \epsilon, \epsilon' \in \{ \pm 1\}.$

 \noindent By a ping-pong argument it follows that $G$ is discrete and free over the generating set  $S$; moreover, its Dirichlet domain with respect to $\tilde{o}$ is
 \vspace{-3mm}

 $$D (G, \tilde{o}) = 
  \Hy^2 \setminus \bigcup_{g_n \in S} \left( A(g_n , \tilde{o}) \cup A(g_n ^{-1}, \tilde{o}) \right)^o$$

\vspace{-1mm}
\noindent  If the axes of the hyperbolic generators do not intersect and the domains $A(g_n^{\pm 1}, \tilde{o}) $ accumulate to one boundary point $\zeta$ (or to different boundary points $E=\{\zeta_k \}$, all defining the same end of the quotient $X=G \backslash \Hy^2$) then the resulting  surface $X=G \backslash \Hy^2$ is a hyperbolic flute: it has a cusp for every parabolic generator, a funnel for every hyperbolic generator, and an {\em infinite end} corresponding to $\zeta$ (or to the set $E$). For the construction of Schottky groups we will repeatedly make use of the following (cp.  Appendix \ref{apphyp} for a proof):
\vspace{-2mm}

\begin{lemma}
\label{lemmahyp1}
Let $\tilde{o} \in \Hy^2$, and let $C, C'$ two ultraparallel geodesics (i.e. with no common point in $\Hy^2 \cup \partial \Hy^2$) such that $d(\tilde{o} , C)=d(\tilde{o} , C')$. Then:

\noindent (i) there exists a unique hyperbolic isometry $g$ with axis perpendicular to $C,C'$ and
 such that $g(C)=C'$;

\noindent (ii) $g^{-1}\tilde{o}$ and $g\tilde{o}$ are obtained, respectively, by the hyperbolic reflections of $\tilde{o}$ with respect to $C,C'$;

\noindent (iii) the Dirichlet domain $D(g,\tilde{o})$ has boundary $C \cup C'$.
\end{lemma}

\begin{example}
\label{exasymmflute}
{\em  The Asymmetric Hyperbolic Flute}

\noindent We construct a hyperbolic flute $X=G\backslash \Hy^2$ with two rays $\alpha, \alpha' $ having same origin such that:

\noindent (a) $\alpha' \prec_G \alpha \not \prec_G \alpha'$
(i.e. $\alpha' \prec \alpha \not \prec \alpha'$);
therefore, $\alpha \not\approx_G \alpha'$ and  $B_{\alpha} \neq B_{\alpha'}$;

\noindent (b) $d_{\infty} (\alpha, \alpha') = \infty$.
\end{example}

\begin{figure}[h]
\label{figflute}
\hspace{3mm}
\includegraphics[width=120mm, height=65mm]{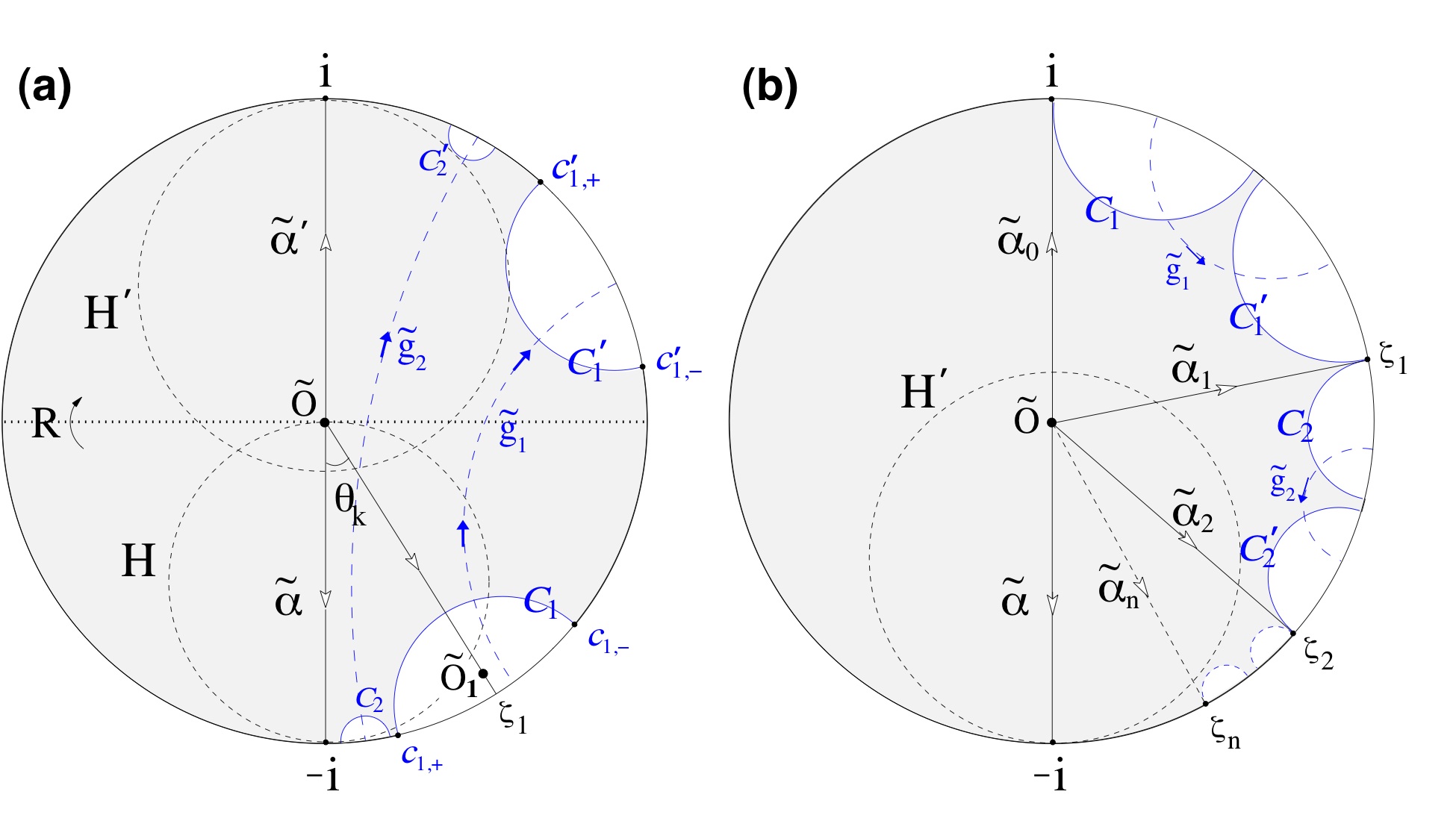}
\vspace{-4mm}
\caption{{\em Asymmetric and Twisted Flutes}}
\end{figure}

\vspace{-1mm}
\noindent  We use the disk model for $\Hy^2$ with origin $\tilde{o}$. Let $\tilde{o}'=-\frac{i}{10}$, and
consider the geodesics  $\widetilde{\alpha}=[\tilde{o}, -i]$,  $\widetilde{\alpha}'=[\tilde{o}, i]$.
Then, let $R$ be  the reflection with respect to the real axis, 
 and consider the horoballs $H=H_{\tilde{\alpha}^+} (\tilde{o})$ and $H'=H_{\tilde{\alpha}^{\prime +}} (\tilde{o}') \supset R(H)$; finally,  choose some positive sequence $\epsilon_k \searrow 0$. \\
Let $[\tilde{o},\zeta_1]$ be a ray making angle $\vartheta_1$ with $\widetilde{\alpha}$,  let $\tilde{o}_1$ be the point on  $[\tilde{o},\zeta_1]$ such that $d(\tilde{o}_1, H)=\epsilon_1$, and let $C_1$ be the hyperbolic perpendicular bisector of the segment $[\tilde{o},\tilde{o}_1]$, with extremities $c_{1,+}$ and $c_{1,-}$, cp.  Figure \ref{figflute}.a. Notice that, as $\epsilon_1 >0$ the circle $C_1$ does not intersect $\widetilde{\alpha}$ (the extremity $c_{1,+}$ closest to $\tilde{\alpha}^+$ coincides with $\tilde{\alpha}^+$ if and only if $\tilde{o}_1 \in \partial H$).
Then, consider $R(C_1)$ and rotate it clockwise around $\tilde{o}$ until it is tangent to $H'$: call this new geodesic  $C_1'$ and its extremities $c_{1,+}', c_{1,-}'$.
\vspace{-3mm}

\noindent Let now $g_1$ be the hyperbolic isometry given by Lemma \ref{lemmahyp1}, with axis $\tilde{g}_1$ \linebreak perpendicular to $C_1,C'_1$, and such that $g_1 (C_1)=C'_1$ and $g_1^{-1} (\tilde{o})=\tilde{o}_1$. \linebreak
Then, construct $g_2$ analogously: that is,  choose a ray  $[\tilde{o},\zeta_2]$, for some $\zeta_2$ between $\widetilde{\alpha}^+$ and $c_{1,+}$,  making angle $\vartheta_2 < \vartheta_1$ with $\widetilde{\alpha}$; call $\tilde{o}_2$ the point on  $[\tilde{o},\zeta_2]$ with $d(\tilde{o}_2, H)=\epsilon_2$, and then let $C_2$,  $C'_2$, $\tilde{g}_2$ etc. as before.  
Repeating inductively this construction we obtain the infinitely generated group $G=<g_1, g_2, ..., g_k, ...>$.   
\vspace{1mm}

\noindent Moreover, choosing $\vartheta_{k+1} \ll \vartheta_k$, we can make the following conditions satisfied:
\vspace{-6mm}

\begin{equation} 
\label{eqschottky}
A(g_n^{\epsilon}, \tilde{o}) \cap A(g_m^{\tau}, \tilde{o}) = \emptyset \mbox{ for all } n \neq m 
 \mbox{ and } \epsilon, \tau \in \{ \pm 1\}
\end{equation} 

\vspace{-4mm}

\begin{equation} 
\label{eqdistance}
U_n ( \widetilde{\alpha} \cup \widetilde{\alpha}' ) \cap A(g_n^{\epsilon}, \tilde{o})  = \emptyset \mbox{ for all } n \in \N \mbox{ and }  \epsilon \in \{ \pm 1\}
\end{equation} 

\noindent where $U_n ( \widetilde{\alpha} \cup \widetilde{\alpha}' )$ is the tubular neighbourhood of $\widetilde{\alpha} \cup \widetilde{\alpha}' $ of width $n$. \\
Condition (\ref{eqschottky}) says that  $ G$ is  a discrete Schottky group.
The quotient manifold $X=G \backslash \Hy^2$ is a hyperbolic flute, with infinite end corresponding to the set $E=\{ \widetilde{\alpha}^+, \widetilde{\alpha}^{\prime +} \}$. Let $\alpha$ and $\alpha'$ be projections of  $\widetilde{\alpha}, \widetilde{\alpha}'$ to $X$, with common origin $o$: they are rays, as their lifts stay in $D(G, \tilde{o})$ by construction. 
\vspace{2mm}

 {\em Proof of Properties \ref{exasymmflute}(a)\&(b).}
 
\noindent  We have $\alpha \succ_G \alpha' $ as $g_n \widetilde{\alpha}^+ \rightarrow \widetilde{\alpha}^{\prime+}$ and    
$B_{\tilde{\alpha}} (\tilde{o}, g_n^{-1} \tilde{o}) \rightarrow 0$, by construction.
\linebreak 
On the other hand,  for every sequence $h_k \in G$ such that $h_k \widetilde{\alpha}^{\prime+}\rightarrow \widetilde{\alpha}^+$, the points $h_k^{-1} \tilde{o}$ definitely lie in some of the attractive domains  $A(g_n, \tilde{o})$, which are exterior to $H'$: thus, $B_{\tilde{\alpha}'} (\tilde{o}, h_k^{-1} \tilde{o}) \! \geq \! \frac{1}{10}$ and does not tend to $0$.
This \nolinebreak proves that $\alpha \not \prec_G \alpha'$. The other assertions in (a) follow from the construction of $G$ and Theorem \ref{teoreq}. For (b), assume that $d_{\infty} (\alpha, \alpha') < M$: then we could find arbitrarily large $t, t'$ and $g_t \in G$ such that $d(\widetilde{\alpha} (t), g_t \widetilde{\alpha}' (t')) <M$. Let then $g_{n(t)}$ be the generator such that $g_t \widetilde{\alpha}' \subset A (g_{n(t)}^{\epsilon}, \tilde{o})$, for some $\epsilon \in \{\pm1\}$. \linebreak
By (\ref{eqdistance}) we deduce that 
$d(\widetilde{\alpha}(t), g_t \widetilde{\alpha}' (t')) \geq d(\widetilde{\alpha}, A(g_{n(t)}^{\epsilon}, \tilde{o}) ) \geq n(t)$
which shows that  we necessarily have $n(t)=n$ for infinitely many, arbitrarily large $t$. Hence 
\vspace{-3mm}

$$\limsup_{t \rightarrow +\infty} d(\widetilde{\alpha} (t), g_t \widetilde{\alpha}' (t')) \geq 
\limsup_{t \rightarrow +\infty}  d(\widetilde{\alpha} (t), A(g_{n}^{\epsilon}, \tilde{o}) ) =\infty$$

\vspace{-1mm}
\noindent a contradiction.\qed

\begin{example}
\label{exsymmflute}
{\em  The Symmetric Hyperbolic Flute}

\noindent We construct a hyperbolic flute $X=\widehat{G}\backslash \Hy^2$ with two rays $\alpha, \alpha' $  having same origin such that:

\noindent (a) $\alpha \prec_{\widehat{G}}\succ \alpha'$
(i.e. $\alpha \prec \succ  \alpha'$);
therefore,  $B_{\alpha} = B_{\alpha'}$;

\noindent (b) $\alpha \not\approx_G \alpha'$;

\noindent (c) $d_{\infty} (\alpha, \alpha') = \infty$.
\end{example}

\vspace{-1mm}
\noindent Let $G=<g_1,...,g_n,...>$ be the group constructed in the Example \ref{exasymmflute}, and  let $S$ be the symmetry with respect to $\tilde{o}$.
Then, for every $n$, consider the hyperbolic translation $\hat{g}_n$ having axis $S[\tilde{g}_n]$ and attractive/repulsive domains $A(\hat{g}_n^{\pm1}, \tilde{o}) = S[A(g_n^{\pm1}, \tilde{o}) ]$, and define $\hat{G}=<g_1,\hat{g}_1,...,g_n,\hat{g}_n,...>$.

\noindent Notice that, by symmetry, all these generators again satisfy the conditions (\ref{eqschottky}) and (\ref{eqdistance}), so $\hat{G}$ is a discrete Schottky group. Again, the quotient manifold $X=\hat{G} \backslash \Hy^2$ is a hyperbolic flute, with infinite end corresponding to the set $E=\{\widetilde{\alpha}^+, \widetilde{\alpha}^{\prime+} \}$  and, with the same notations as above, the projections $\alpha$ and $\alpha'$ on $X$ are rays.
\vspace{3mm}

  {\em Proof of Properties \ref{exsymmflute}(a),(b)\&(c).}
  
\noindent
We deduce as before that $ \alpha \succ_{\hat{G}} \alpha' $; but now  we also have the sequence $\hat{g}_n$ such that $\hat{g}_n \widetilde{\alpha}^{\prime+}  \rightarrow \widetilde{\alpha}^+$ and $B_{\alpha'} (\tilde{o}, \hat{g}_n^{-1} \tilde{o}) \rightarrow 0$; so $\alpha' \succ_{\hat{G}} \alpha$ too.
As the rays $\alpha$ and $\alpha'$ have a common origin,  Theorem \ref{teoreq} implies that $B_{\alpha} = B_{\alpha'}$. Again assertion (b) follows by construction, and (c) is proved as before.\qed
\vspace{5mm}

\begin{example}
\label{extwistedflute}
{\em The Twisted Hyperbolic Flute}

\noindent We construct a hyperbolic flute $X=G\backslash \Hy^2$ with a family of rays $\alpha_n$ having same origin and converging to a ray $\alpha$ such that:

\noindent (a) $\; \alpha_n \approx_G \alpha_m$ $ \forall n,m$; 
 therefore, $d_{\infty} (\alpha_n,\alpha_m) <  \infty$ and $B_{\alpha_n} = B_{\alpha_m}$ $ \forall n,m$;

\noindent (b) $\; d_{\infty} (\alpha_n,\alpha) =\infty$ $\forall n$; 

\noindent (c) $B_{\alpha_0} = \lim_{n\rightarrow +\infty }B_{\alpha_n}  \neq B_{\alpha}$.
\end{example}
 
\vspace{-1mm}
\noindent Again, in the disk model for $\Hy^2$ with origin $\tilde{o}$,  consider a sequence of boundary points $\zeta_0= i$, $\zeta_n = e^{i \vartheta_n}$, for a decreasing sequence $\frac{\pi}{2} \geq \vartheta_n \searrow -\frac{\pi}{2}$. \linebreak 
Then, for every $n \geq 1$ choose a pair of ultraparallel geodesics   $C_n, C'_n$ such that \linebreak  $d(\tilde{o}, C_n)= d(\tilde{o}, C'_n)=d_n$,  each cointained in the disk sector $[\zeta_{n-1},\tilde{o}, \zeta_{n}]$, and  with  points at infinity respectively equal to $\zeta_{n-1}$,  $\zeta_{n}$. Finally, let $g_n$ be the hyperbolic isometry with $g_n (C_n)=C'_n$ whose axis is perpendicular to $C_n,C'_n$, given by Lemma \ref{lemmahyp1}, cp. Figure \ref{figflute}.b, and set $\widetilde{\alpha}_n = [\tilde{o}, \zeta_n]$, $\widetilde{\alpha}= [\tilde{o}, -i]$. \\
Moreover, if $H' = H_{\tilde{\alpha}^+} (\tilde{o}')$ for $\tilde{o}'=\frac{i}{10}$, we can choose the  $d_n \gg0$ in order that the following condition is satisfied:
\vspace{-3mm}

\begin{equation}
\label{eqnonintersecting}
H'  \cap A(g_n^{\pm 1}, \tilde{o})  = \emptyset \mbox{ for all } n
\end{equation}

\vspace{-1mm}
\noindent Define $G$ as the group generated by all the $g_n$. Again, $G$ is an infinitely generated Schottky group, and the quotient manifold $X=G\backslash \Hy^2$ is a flute whose infinite end corresponds to  the set $E=\{ \widetilde{\alpha}^{+}, \widetilde{\alpha}_n^+ \; | \; n \geq 0 \}$. The projections $\alpha_n$ and $\alpha$ of all the $\widetilde{\alpha}_n, \widetilde{\alpha}$ on $X$ are rays, by construction, such that $\alpha_n \rightarrow \alpha$. 
\vspace{3mm}
 
 {\em Proof of Properties \ref{extwistedflute}(a),(b)\&(c).}
  
\noindent The rays $\alpha_n$ are all $G$-equivalent by construction, as $ \widetilde{\alpha}_n^+=g_n \widetilde{\alpha}_{n-1}^+$ for all  $n$. \linebreak The other assertions in (a) follow from the discussion after Definition \ref{defieqG}\linebreak
(actually, as we are in strictly negative curvature,  we have $d_{\infty}(\alpha_n, \alpha_m) =0$).\linebreak
On the other hand, by (\ref{eqnonintersecting}),  all the images by $G$ of $\widetilde{\alpha}_n$ are exterior to the horoball $H'$, exceptly for $\widetilde{\alpha}_n$ itself; thus if $s\gg0$  we have $d(g\widetilde{\alpha}_n, \widetilde{\alpha} (s)) \! > \!s$ for all $g$. It follows that $d_{\infty} (\alpha_n, \alpha) \geq \frac12 \limsup_{s \rightarrow +\infty} \inf_{g\in G} d(g\widetilde{\alpha}_n, \widetilde{\alpha}(s)) =+\infty$. \linebreak
To conclude we have to prove that $B_{\alpha_0} \neq B_{\alpha}$, and by Theorem \ref{teoreq} it is enough to show that 
$ \alpha \not\succ_G \alpha_0$. But for any sequence $h_n$ with $h_n \widetilde{\alpha}^+ \rightarrow \widetilde{\alpha}_0^+$  we have $B_{\tilde{\alpha}}(\tilde{o}, h_n \tilde{o})< -\frac{1}{10}$, since by construction this is true for all nontrivial $g$ in $G$.\qed

\begin{remark}
\label{remcontinuity}
{\em 
The discontinuity (c) can be interpreted geometrically as follows. \\
Consider the maximal horoballs  $H_{\tilde{\alpha}^+}^{max} (o'), H_{\tilde{\alpha}_n^+}^{max} (o')$, for the projection   $o'$  of \nolinebreak  $\tilde{o}'$.
It is easy to see that $H_{\tilde{\alpha}^+}^{max} (o')= H_{\tilde{\alpha}^+} (\tilde{o}')$, as all the  $g\tilde{o}'$, for $g \neq 1$, stay far away from $H'$, by construction.
Moreover, since $o' \in \alpha_0$ and $\alpha_0$ is a ray, we  also deduce that $H_{\tilde{\alpha}_0^+}^{max} (o') = H_{\tilde{\alpha}_0^+} (\tilde{o}')$ precisely.
Now   $B_{\alpha_n}(o,o') = B_{\alpha_0}(o,o')$, so formula \nolinebreak   (\ref{eqbusemannb})  shows that 
$d(\tilde{o}, H^{max}_{\tilde{\alpha}_n^+} (o')) = d(\tilde{o}, H^{max}_{\tilde{\alpha}_0^+} (o'))$; then,  by rotational symmetry,  $H_{\tilde{\alpha}_n^+}^{max} (o')$ is  the horoball centred at $\tilde{\alpha}_n^+$ having the same Euclidean radius as $H_{\tilde{\alpha}_0^+} (\tilde{o}')$.
Therefore the discontinuity can be read in terms of a discontinuity in the limit of the maximal horoballs: in fact, the $H_{\tilde{\alpha}_n}^{max} (o')$'s converge  for $n\rightarrow \infty$  to $H_{\tilde{\alpha}^+} (-\tilde{o}')$, which is strictly smaller than the maximal horoball  $H_{\tilde{\alpha}^+} (\tilde{o}')$ of the limit ray.}
\end{remark}
\vspace{3mm}

\begin{example}
\label{exladder}
{\em The Hyperbolic Ladder}

\noindent We construct a hyperbolic ladder which is a Galois covering $X \rightarrow \Sigma_2$ of a hyperbolic surface of genus $2$,  with automorphisms  group  $\Gamma \cong \mathbb{Z}$, such that:

\noindent (a)  $X$ has distance-asymptotic rays $\alpha, \alpha'$ with $\alpha \prec\succ \alpha'$, but $B_{\alpha} \neq B_{\alpha'}$;

\noindent (b)  ${\cal B}X$ consists of 4 points;

\noindent (c) $\partial X$ consists of a continuum of points;

\noindent (d)  the limit set  $L\Gamma= \overline{\Gamma x_0} \cap \partial X$  depends on the choice of the base point $x_0$, and for some $x_0$ it is included in  $\partial X - {\cal B}X$.
\end{example}

\noindent We construct $X$ by glueing infinitely many pairs of hyperbolic pants.\\
The following properties of hyperbolic pants are well-known:

\begin{lemma}[\cite{fathi}, \cite{thurston}] 
\label{lemmapant}
Let $H^+, H^-$ be two identical  right-angled hyperbolic hexagons with alternating edges labelled respectively by $a^{\pm}, b^{\pm}, c^{\pm}$ and opposite edges $\alpha^{\pm} ,  \beta^{\pm},  \gamma^{\pm}$.
Let $P$ the hyperbolic pant obtained by glueing them along $a^{\pm}, b^{\pm}, c^{\pm}$;  the identified edges $a,b,c$ are called the  {\em seams} of $P$, and the resulting boundaries $\alpha= \alpha^+\cup \alpha^-,  \beta= \beta^+\cup \beta^-,  \gamma=\gamma^+ \cup \gamma^-$ of $P$ are closed geodesics  called the  {\em cuffs}. The seams are the shortest geodesic segments connecting the cuffs of $P$ and, reciprocally, the cuffs are the shortest  ones connecting the seams.

\end{lemma}

\begin{figure}[h]
\label{figladder}
\hspace{-10mm}
\includegraphics[width=14cm, height=5cm]{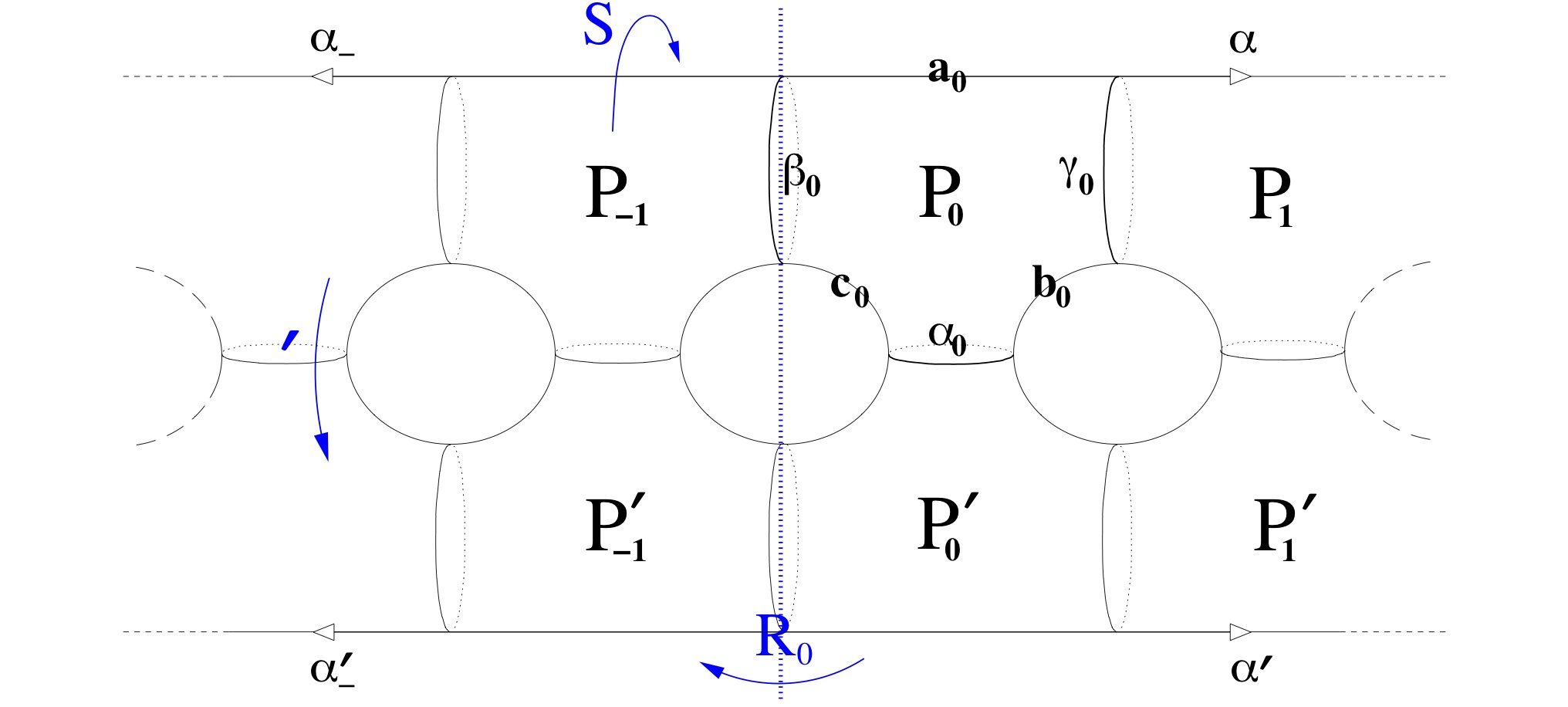}
\caption{The Hyperbolic Ladder}
\end{figure}

\noindent Now, we start from infinitely many copies $P_n$, $P'_n$, for $n \in \mathbb{Z}$, of the same pair of pants $P$, and we assume that $\ell (b) = \ell (c) =L > \ell = \ell (a)$.
 We glue them as in figure \ref{figladder}, by identifying via the identity the cuffs $\alpha_n$ with $\alpha'_n$, and the cuffs $\beta_n$, $\beta'_n$ with $\gamma_{n-1}$, $\gamma'_{n-1}$ respectively (with no twist), obtaining a complete hyperbolic surface  $X = N \backslash \Hy^2$.
Remark that, if $\Sigma_2 =G \backslash \Hy^2$ is the hyperbolic surface obtained from $P_0 \cup P'_0$ by identifying $\alpha_0$ to $\alpha'_0$ and $\beta_{0}$, $\beta'_0$ respectively to $\gamma_0$, $\gamma_0'$, there is a natural covering projection  $X \rightarrow  \Sigma_2$,  with automorphism group \nolinebreak  $\Gamma \cong \mathbb{Z} \cong G/N$. The group $\Gamma$ acts on $X$ by ``translations'' $T_k$, sending $P_n \cup P'_n$ into $P_{n+k} \cup P'_{n+k}$.

\noindent We define 
$\alpha     \!=\!  \bigcup_{k \geq 0} a_k$,
$\alpha_- \!= \! \bigcup_{k < 0} a_k$, 
$\alpha'    \!=\!  \bigcup_{k \geq 0} a'_k$,
$\alpha'_- \!= \! \bigcup_{k < 0} a'_k$
and set 
$\mathcal{A}=   \alpha  \cup \alpha_-$, $\mathcal{A}'=\alpha'    \cup \alpha' _- $.
 Notice that the surface $X$ is also endowed of:
\vspace{1mm}

\noindent --  a natural {\em flip symmetry}, denoted $'$,  obtained by sending a point in $P_k$ to the corresponding point in $P'_k$; let $\mathcal{F}= Fix(')$ and  call the  {\em top} and the {\em bottom} of $X$ the (closure of) the two connected components of $X - \mathcal{F}$ interchanged by $'$;
\vspace{1mm}

\noindent -- a natural {\em mirror symmetry} $S$, obtained interchanging each point on a pant $P_k$ (resp. $P'_k$) with the corresponding point lying on the same pant, but on the opposite hexagon; 
if $\mathcal{M}= \bigcup_{k\in \mathbb{Z}} b_k \cup b'_k \cup c_{k} \cup  c'_{k}$, we have $Fix(S)= \mathcal{A} \cup \mathcal{A}' \cup \mathcal{M}$, and 
we will call the {\em back} and the {\em front} of $X$ the closure of the two connected components of $X - Fix(S)$ interchanged by $S$;
\vspace{1mm}

\noindent -- a group of {\em reflections} $R_n$ with respect to $\beta_n \cup \beta'_n$, exchanging $P_{n+k} \cup P'_{n+k} $ with  $P_{n-k-1} \cup P'_{n+k-1}$.

\begin{lemma}
\label{propladder}
${}$

\noindent (i) Every  minimizing geodesic does not cross twice  neither $\mathcal{A}$, $\mathcal{A}'$, $\mathcal{F}$ nor $\mathcal{M}$;

\noindent (ii) every quasi-ray is strongly asymptotic to one of the four rays $\alpha, \alpha_-,\alpha', \alpha'_-$.

\end{lemma}

{\em Proof.}

\noindent (i) Assume that $\gamma$ is a minimizing geodesic between $x$ and $y$, crossing $\mathcal{A}$ twice,  at two points $x_1,y_1$. Break it as  $\gamma =   \gamma_1 \cup [x_1,y_1] \cup \gamma_2 $, where  $[x_1,y_1] $ is the subsegment between $x_1$ and $y_1$. Then, using the mirror symmetry $S$, we would obtain a curve $\hat{\gamma} =   \gamma_1 \cup S[x_1,y_1] \cup \gamma_2 $ of same length, still connecting $x$ to $y$, but singular at $x_1$ and $y_1$; hence, it could be shortened, which is a contradiction. The proof is the same for $\mathcal{A}', \mathcal{B}$, and using the flip simmetry $'$ one analogusly proves that  a minimizing geodesic cannot cross twice  $\mathcal{F}$. \\
For (ii), let us first show that, if $\gamma$ is a quasi-ray included, say, in the top-front of $X$, then either $d_{\infty} (\gamma, \alpha)=0$ or $d_{\infty} (\gamma, \alpha_-)=0$.
Actually,
assume that $p_n = \gamma(t_n)$ is a sequence such that $d(p_n, \mathcal{A}) > \epsilon$, for  $n \geq 0$ and  $t_n \rightarrow \infty$. Consider the projections $q_n$ of $p_n$ on $\mathcal{A}$, which we may assume at distance $d(q_n, q_{n+1}) \gg 0$;
 as $\gamma$ is included in a simply connected open set of $X$ containing the bi-infinite geodesic $\mathcal{A}$, we can  use hyperbolic trigonometry  (cp. Lemma \ref{lemmalength} in the \S\ref{apphyp}) to deduce that $\ell (\gamma |_{[t_n, t_{n+1}]}) \geq  q_n q_{n+1}  + \delta (\epsilon)$, for a  universal function $\delta (\epsilon)>0$. \\
 As $p_0 p_n \leq q_0 q_n + 2 \mbox{diam}(P)$,  we obtain
\vspace{-6mm}

$$ \Delta(\gamma|_{[t_0, t_{N}]})
\geq \sum_{n=0}^{N-1}  q_n q_{n+1}  + N \delta (\epsilon) -  q_0 q_N  - 2 \mbox{diam}(P)
=N \delta (\epsilon) - 2 \mbox{diam}(P)$$

\vspace{-2mm}
\noindent which diverges as $N \rightarrow \infty$; so, $\Delta (\gamma)$ is not bounded, a contradiction. As $\epsilon$ is arbitrary, this shows that  $\gamma$ is strongly asymptotic either to $\alpha$ or to $\alpha_-$.
Finally, if $\gamma$ is a quasi-ray which is not included in the top-front of $X$, we can use the symmetries $S$ and $'$ to define, from $\gamma$, a curve $\hat{\gamma}$ fully included in the top-front of $X$, by reflecting the subsegments which do not lie in the top-front of $X$. This new curve still has bounded excess (as it has the same lenght as $\gamma$ on every interval, and the distance between endpoints reduces at most of $2 \mbox{diam}(P)$) so, as we just proved, it is strongly asymptotic  either to $\alpha$ or to $\alpha_-$. In particular, $\hat{\gamma}$ finally does not intersect $C$; so, $\gamma|_{[t_0, +\infty]}$,  for some $t_0 \gg 0$, is   included in an $\epsilon$-neighborhood of $\mathcal{A}$, for arbitrary $\epsilon$, and therefore it is strongly asymptotic to one of the four rays $\alpha, \alpha_-, \alpha', \alpha'_-$.\qed
\vspace{4mm}

{\em Proof of \ref{exladder}(a),(b),(c)\&(d) of Example \ref{exladder}}.

\noindent {\em (a)}
The geodesic segments $a_n$ are the shortest curves connecting the cuffs $\beta_n, \gamma_n$ of $P_n$: this implies that $\alpha$  cannot be shortened, so it is a ray; similarly for  $\alpha'$. \linebreak
Let  now $x_0 = a_0 \cap  \beta_0$, $x_n=T_n(x_0)$ and let $x'_n$ be their flips; finally, consider a sequence of minimizing segments $\eta_n = [x_0, x'_{2n}]$ and their inverse paths $-\eta_n$. \linebreak 
By (i) above, we know that $\eta_n$ is  included in the front (or the back) of $X$; moreover, it can be broken as $\eta_n = \eta^t_n \cup \eta^b_n$ where  $\eta^t_n, \eta^b_n$ respectively are subsegments in the top and in the bottom of $X$ meeting at some $p_n \in \mathcal{F}$ . Therefore, each of these segments stays in a simply connected  open set of $X$, isometric to an open set of $\Hy^2$; then,  since $d(p_n, \alpha )= d(p_n, \alpha')  < \mbox{diam}(P)$ we can apply standard hyperbolic trigonometry to 
deduce that  $\eta_n$ makes an angle $\vartheta_n$, with either $\alpha$ or  $\alpha'$, such that
\vspace{-5mm}

$$\tan \vartheta_n \leq \frac{\tanh (diam(P))}{\sinh (n \ell)} \rightarrow 0, \; \mbox{ for } n \rightarrow \infty.$$

\vspace{-1mm}
\noindent By possibly replacing $\eta_n$ with  $R_{n} (-\eta_{n})'$, we find  a sequence of minimizing  segments $[x_0, x'_{2n}] \rightarrow \alpha$, hence $\alpha' \succ \alpha$. The converse relation $\alpha \succ \alpha'$ is analogous. Let us now show that $B_{\alpha} \neq B_{\alpha'}$. It is enough to show that $B_{\alpha} (x_0,x'_0)>0$; then clearly, by the flip symmetry, we will deduce  $B_{\alpha'} (x_0,x'_0) = B_{\alpha} (x'_0,x_0) <0$. 
\linebreak
Let us compute $B_{\alpha} (x_0,x_0') = \lim_{n \rightarrow \infty}  x_0 x_n  - x_n x_0' $. 
Let $\nu_n=[x_n,x_0']$ be a minimizing segment intersecting $\mathcal{F}$ at some   $p \in\hat{\alpha}_k$, and break it as
$\nu_n=\nu^t_n \cup \hat{\nu}_n \cup \nu^b_n$
where $\hat{\nu}_n $ is the maximal subsegment of $\nu_n$ included in $P_k \cup P'_k$; then, 
\vspace{-3mm}

$$ x_n x_0'  \geq \ell(\nu^t_n) +d(\gamma_k, \beta'_k) + \ell( \nu^b_n) \geq   (n-1)\ell + 2L \geq (n+1) \ell$$

\vspace{-1mm}
\noindent while  clearly $x_0 x_n = n\ell$; so, $B_{\alpha} (x_0,x_0') \geq \ell$.
\vspace{1mm}

\noindent {\em (b)} One proves analogously that $\alpha_-$ and $ \alpha'_-$ are rays defining different Busemann functions, while it is clear that $B_{\alpha}$ and $B_{\alpha'}$ are different  from $B_{\alpha_-}, B_{\alpha'_-}$. Therefore  ${\cal B}X$ has at least 4 points. On the other hand, by   Proposition   \ref{propladder}(ii), every quasi-ray in $X$ is strongly asymptotic to one of the four above, thus defining the same Busemann function. This shows that ${\cal B}X$  has precisely four points. \\
 
\vspace{-3mm}
\noindent {\em (c)-(d)} Clearly, the orbits $\Gamma x_0$ and $\Gamma x'_0$ accumulate to $B_{\alpha}$ and $B_{\alpha'}$.
Let now $x(t)$ be a continuous curve from $x_0=x(0)$ to $x'_0=x(1)$, and set  $x_n(t)=T_n (x(t))$. For any fixed $t$,  let $B_{(x_n)(t)}$ be the limit of (a subsequence of) $x_n(t)$, for $n\rightarrow \infty$. The family $B_{(x_n)(t)}$ defines a continuous curve in $\partial X$ connecting $B_{\alpha}$ to $B_{\alpha'}$, as
$\|B_{(x_n)(t)} - B_{(x_n)(s)}\|_{\infty} \leq 2 d(x_n(t), x_n(s))$; since it is non-constant,   its image is an uncountable subset of $\partial X$.
It remains to exhibit an orbit accumulating to a point of $\partial X \setminus {\cal B}X$. Let   $y_0\in  \alpha_0$: we affirm that $y_n=T_n  y_0$ is such an orbit. Actually, if $y_n$ converged to one of the four Busemann functions, say $B_{\alpha}$, then we would also have $y_n =y'_n  \rightarrow B_{\alpha'}$, as the flip symmetry preserves the orbit and exchanges $\alpha$ with $\alpha'$. Hence we would get $B_{\alpha}=B_{\alpha'}$, a contradiction.\qed

\begin{remark}
{\em The surface $X$ is quasi-isometric to $\mathbb{Z}$, hence it is a Gromov-hyperbolic metric space. Its {\em boundary as a Gromov-hyperbolic space}  $X^g (\infty)$ (cp. \cite{bridson}, \cite{papa}) consists of two points. So, the  Busemann boundary and the horoboundary prove to be  finer invariants than $X^g (\infty)$ (as they are not defined up to bounded functions, so they are not invariant by quasi-isometries).
}
\end{remark}

\appendix
\section{Appendix}

\subsection{Rays on Riemannian manifolds.}
\label{appcorays}

\begin{lemma}
\label{lemmacoray}
Let $\beta$ be a quasi-ray and $x,y\in X$ such that $B_{\beta} (x,y) \!= \!d(x,y).\!$ \nolinebreak
Then:

\noindent (i) $x$ and $y$ minimize the distance between the horospheres  $H_{\beta^+}(x)$ and  $H_{\beta^+}(y)$;

\noindent (ii) $y$ is the {\em only} projection to $H_{\beta^+} (y)$ of every   $z \!\in\! [x,y]$, exceptly possibly for 
\nolinebreak $x$.
\end{lemma}

\vspace{-1mm}
\noindent {\bf Proposition \ref{propcoray}}
{\em
For any quasi-ray  $\beta$ we have: 
 $\;B_{\beta} (x, y) = d(x, y)\; \Leftrightarrow \; \overrightarrow{xy}\prec \beta$. \linebreak
 In particular, if $B_{\beta} (x, y) = d(x, y)$, the extension of any minimizing segment $[x,y]$ beyond $y$ is always a ray.
 }

\vspace{2mm}  
\noindent {\bf Theorem \ref{teorcoray}} \hspace{1mm}
{\em 
Assume that  $\alpha, \beta$ are rays in $X$ with origins $a,b$ respectively.  \linebreak
 The following conditions are equivalent: 
\vspace{1mm}

\noindent (a) $B_{\alpha} (x,y) \! = \! B_{\beta} (x,y)$  $ \forall x, y \in X$;
\vspace{1mm}

\noindent (b)  $\alpha \prec\succ \beta$ and $B_{\alpha} (a,b) = B_{\beta} (a,b)$;
\vspace{1mm}

\noindent (c)  $\alpha$ and $\beta$ are visually equivalent from every  $o \in X$.
}

\vspace{3mm}
{\em Proof of Lemma \ref{lemmacoray}}. 
(i)  follows from the fact that any two points $x',y'$ respectively in $H_{\beta^+}(x), H_{\beta^+}(y)$ satisfy
$d(x',y') \geq B_{\beta} (x',y')=B_{\beta} (x,y) =d(x,y)$. 
In particular, $y$ is a projection to $H_{\beta^+}(y)$ of any point $z \in [x,y]$, as
\vspace{-3mm}

$$ xz + zy  = xy  = d(x, H_{\beta^+}(y)) \leq xz + d(z, H_{\beta^+}(y)). $$

\vspace{-1mm}
\noindent Moreover,  let $z \in [x,y]$, $z \neq x$, and assume that $q$ is a projection  of $z$ on $H_{\beta} (y)$ different from $y$. Then, the angle between $[x,z]$ and $[z,q]$ would be different from $\pi$; hence $xq < xz +zq $ and 
\vspace{-3mm}

$$d (x, H_{\beta^+}(y)) < xz + zq =  xz + zy  = d(x,  H_{\beta^+} (y))$$

\vspace{-1mm}
\noindent a contradiction.\qed

\vspace{4mm}
\emph{Proof of Proposition \ref{propcoray}.}
\hspace{-1mm}
Let $\alpha \!=\!  \overrightarrow{xy}$, with  $x\!=\!\alpha(0)$, $y\!=\!\alpha (\overline{s})$. 
Assume \nolinebreak $\alpha \!\prec\! \beta$.
So, there exist minimizing geodesic segments $\alpha_n=[a_n, b_n] \rightarrow \alpha$ such that $a_n \!=\!\alpha_n(0) \rightarrow x$, $b_n\!=\!\alpha_n(s_n)\!=\! \beta (t_n)  \rightarrow \beta^+$, for  sequences $s_n,t_n \rightarrow +\infty$. \linebreak 
 Let $s$ be fixed and $\epsilon$  arbitrary. There exists $N(s,\epsilon)$ such that $d(\alpha_n(s), \alpha(s)) < \epsilon$ and $d(a_n, x)< \epsilon$ for $n > N(s, \epsilon)$; therefore
 \vspace{-3mm}
 
$$B_{\beta} (x, \alpha(s)) = \lim_{n \rightarrow \infty}  x b_n  -  b_n  \alpha(s) 
 \eqsim_{\epsilon}   \lim_{n \rightarrow \infty}  x  b_n -  b_n  \alpha_n(s) =s $$

 \vspace{-1mm}
\noindent and as $\epsilon$ is arbitrary, this shows that $B_{\beta} (x, \alpha(s))  =s= d (x, \alpha(s))$ for all $s$,  hence  $B_{\beta} (x, y)  =d (x, y)$. 
Conversely, assume that $B_{\beta} (x, y)  =d (x, y)$. Then:
 \vspace{-6mm}
 
$$ s = \overline{s} - (\overline{s}-s) \leq B_{\beta} (x, y) -  B_{\beta} (\alpha(s),y) = B_{\beta} (x, \alpha(s)) \leq s   \;\;\; \forall s \in [0,\overline{s}]$$

 \vspace{-2mm}
 \noindent and we deduce that $B_{\beta} (x, x')  =d (x, x')$ for all $x,x'$ on $\alpha$ between $x$ and $y$. \linebreak
Now, fix $0<\epsilon<\overline{s}$  and consider minimizing geodesic segments $\alpha^{\epsilon}_n = [ \alpha (\epsilon) , \beta (n)]$; up to  a subsequence, they converge, for $n \rightarrow \infty$, to a ray $\alpha^{\epsilon}$ which is, by definition, a coray of $\beta$.
So (as we previously proved)
\vspace{-3mm}
 
$$B_{\beta} (\alpha (\epsilon) ,  \alpha^{\epsilon} (s))= 
  B_{\alpha^{\epsilon}} (\alpha (\epsilon) ,  \alpha^{\epsilon} (s)) =s
  \;\;\; \forall s >0$$
  
\vspace{-1mm}
\noindent  But then,  for $s > \epsilon$,  $\alpha (s)$ and $\alpha^{\epsilon} (s-\epsilon)$  are both projections of $\alpha(\epsilon)$ to the horosphere  $H_{\beta^+}(\alpha (s))$  and, by Lemma \ref{lemmacoray}(ii), we know that they coincide. 
This shows that $\alpha^{\epsilon} = \alpha |_{\epsilon, +\infty}$ and that $\alpha^{\epsilon \, \prime}_n  (0)$ tend  to $\alpha' (\epsilon)$, for every fixed $\epsilon>0$; by a diagonal argument we then build a sequence of minimizing geodesic segments $\alpha_k = \alpha^{\epsilon_k}_{n_k}$, for $\epsilon_k\rightarrow 0$ and $n_k \rightarrow +\infty$, such that $\alpha_k  \rightarrow \alpha $.  Thus $\alpha \prec \beta$.\qed

\vspace{4mm}
\emph{Proof of Theorem \ref{teorcoray}.}\\
\noindent Let us show that (a) $\Rightarrow$ (b). Assume that $B_{\alpha}=B_{\beta}$, and let $b=\beta(0), y=\beta(t)$. As $B_{\alpha} (b,y)=B_{\beta} (b,y) = d(b,y)$, we deduce by Proposition \ref{propcoray} that $\beta \prec \alpha$. One proves that  $\alpha \prec \beta$  analogously. 
 
 \noindent Conversely, let us show that (b) $\Rightarrow$ (a).  Assume that $\alpha \prec \beta$, so we have geodesic segments    $\alpha_n=[a_n, b_n] \rightarrow \alpha$ with $a_n=\alpha_n(0) \rightarrow a$, $b_n= \beta (t_n) =\alpha_n(s_n) \rightarrow \beta^+$; moreover,  let  as before $N(s,\epsilon)$ such that $d(\alpha_n(s), \alpha(s)) < \epsilon$ for $n > N(s, \epsilon)$.  Then, for every $x$ and $n>N(s, \epsilon)$:
\vspace{-4mm}

$$  a \alpha(s)  -  \alpha(s) x \; \eqsim_{\epsilon} \;
 s -  \alpha_n(s) x    \;  \leq \;
 s -    \left(   b_n x  -   b_n  \alpha_n(s) \right)  $$

\vspace{-1mm}
\noindent and, as  $b_n \alpha_n(s) = s_n -s$ we deduce that 
\vspace{-6mm}

$$  a \alpha(s)  -  \alpha(s) x   \;\apprle_{\epsilon}\;
    s_n -  b_n x  \; = \;
   (s_n \!-\! t_n)  + ( t_n \! -  b_n x  ) \;\leq\;
    B_{\beta} (a_n,b) + B_{\beta}  (b,x)$$

\vspace{-2mm}
\noindent by monotonicity of the Busemann cocycle. Taking limits for $s \rightarrow \infty$ we deduce that   $B_{\alpha} (a,x) \apprle_{\epsilon} B_{\beta} (a,x) $  for all $x$ and, as $\epsilon$ is arbitrary, $B_{\alpha} (a,x) \leq B_{\beta} (a,x) $.
 From $\beta \prec \alpha$ we deduce analogously that $B_{\beta} (b,x) \leq B_{\alpha} (b,x) $.
Therefore: 
\vspace{-5mm}

$$B_{\beta} (b,x) \leq B_{\alpha} (b,x) =  B_{\alpha} (b,a) + B_{\alpha} (a,x) \leq
B_{\alpha} (b,a) + B_{\beta} (a,b) + B_{\beta}(b,x)$$

\vspace{-1mm}
\noindent and since $B_{\alpha} (b,a) = B_{\beta} (b,a)$ we get the conclusion.

\noindent Let us now prove that (a) $\Rightarrow$ (c). Assume again that $B_{\alpha}=B_{\beta}$, and let $o \in X$.
Let $\gamma$ be a limit of (a subsequence of)   geodesic segments $\gamma_n= [o, \alpha(n)]$; then $\gamma$ is a ray (by the Properties \ref{lemmaexcess}) and, by definition, is a coray to $\alpha$. 
Then,  by Proposition  \ref{propcoray}
\vspace{-5mm}

$$B_{\beta} (o, \gamma(t))=B_{\alpha} (o, \gamma(t)) = d(o, \gamma(t))$$

\vspace{-1mm}
\noindent which, by the same Proposition, also implies that $\gamma \prec \beta$.\\
Finally, let us show that (c) $\Rightarrow$ (a). The functions  $B_{\alpha}(a, \cdot)$ and $B_{\beta}(b, \cdot)$ are Lipschitz, hence differentiable almost everywhere. Let $o$ be a point of differentiability for  both $B_{\alpha}(a, \cdot)$ and $B_{\beta}(b, \cdot)$, and let $\gamma$ be a ray from $o$ which is a coray to $\alpha$ and $\beta$.
Then $B_{\alpha}(o, \gamma(t)) = d(o, \gamma(t))=B_{\beta}(o, \gamma(t))$ for all $t$, which implies that
$grad_{o} B_{\alpha} (a, \cdot)= \gamma' (0) = grad_{o} B_{\beta} (b, \cdot)$. So  $B_{\alpha}(a, \cdot)$ and  $B_{\beta}(b, \cdot)$ are Lipschitz functions whose gradient is equal almost everywhere, hence they differ by a constant and $B_{\alpha}=B_{\beta}$.\qed

\subsection{Rays on Hadamard spaces.}
\label{apphadamard}

\noindent {\bf Proposition \ref{prophadamard}}
{\em Let  $\widetilde{X}$ be  a Hadamard space:
\vspace{1mm}

\noindent (i) if $\alpha, \beta$ are rays, then $d_{\infty} (\alpha, \beta) < \infty \;  \Leftrightarrow \; B_{\alpha}=B_{\beta} \; \;  \Leftrightarrow \; \alpha \prec \beta$.

\noindent  Moreover, two rays with the same origin are Busemann equivalent iff they coincide, so the restriction of the Busemann map 
$B_o \!:\! {\cal R}_o  ( \widetilde{X}) \rightarrow \partial \widetilde{X}$ is injective;

\noindent (ii) for any $o \in \tilde{X}$, the restriction   of the Busemann map $B_o: {\cal R}_o ( \widetilde{X})\rightarrow \partial \widetilde{X}$ is surjective, hence ${\cal B}\tilde{X} = {\cal B}_o \tilde{X} = \partial \tilde{X}$;
\vspace{1mm}

\noindent (iii) the Busemann map $B$ is continuous.
}
\vspace{2mm}

\noindent {\bf Uniform Approximation Lemma \ref{lemmaapprox}}
{\em
Let $\widetilde{X}$ be a Hadamard space. \\
For any compact set $K$ and $\epsilon>0$, there exists a function $T (K, \epsilon)$ such that for any $x,y \in K$ and any ray $\alpha$ issuing from $K$,  we have  $|B_{\alpha} (x,y) - b_{\alpha(t)}(x,y)| \leq \epsilon$, provided that $t \geq T (K, \epsilon)$.
}
\vspace{3mm}

 \emph{Proof of Lemma \ref{lemmaapprox}.} 
 
 \noindent First notice that, by the cocycle condition (holding for $b_{\alpha(t)}$ as well as for $B_{\alpha}$   we can assume that $x=\alpha(0)=a$.
Then,  let $z=\alpha(t)$,   $z'=\alpha(t')$ for $t'>t$, and let us estimate $b_{z'} (a,y) - b_{z} (a,y) = (yz + zz') - yz'$.
Assume that $K \subset B (a,r)$,  denote by  $y'$ the  projection of $y$ on $\alpha$ and consider $\vartheta =  \widehat{y z a}$.
The right triangle $[y,y',z]$ has catheti  $yy' \leq r$ and $zy' \geq t-r$ (as $ay' \leq r$); by comparison with a Euclidean triangle, we deduce that $0<\vartheta \leq \vartheta_0<\pi $ with $\tan \vartheta_0  = r/(t-r)$.
Comparing now the triangle $[y,z,z']$  with an Euclidean triangle $[y_0,z_0,z_0']$ such that $\widehat{y_0z_0z'_0} = \pi - \vartheta_0$ and $y_0z_0=yz$,  $z_0z_0'=zz'$ we deduce that $yz' \geq y_0z_0'$. So, 
\vspace{-2mm}

\begin{equation}
\label{comparison}
 b_{z'} (a,y) - b_{z} (a,y)  = (yz + zz') - yz' \leq  (y_0z_0 + z_0z_0') - y_0z_0' 
\end{equation}

\noindent Now a straightforward computation in the plane shows that this tends to zero uniformly on  $y \in K$,  for  $t \rightarrow \infty$. Actually, consider the projection  $y'_0$ of $y_0$ on the line containing $z_0, z'_0$, and set
$r_0= y_0y'_0$,  $s_0 = z_0z'_0$ and $\rho_0 = y'_0z_0$.
Then, for  $r$ fixed and $t$ tending to infinity, we have
$t+r \geq yz \geq  \rho_0 = yz \cos \vartheta_0 \rightarrow +\infty$
while
$r_0 = \rho_0 \tan \vartheta_0 \leq \frac{r(t+r)}{t-r} $ stays bounded.
Therefore
\vspace{-5mm}

$$(y_0z_0 + z_0z_0') - y_0z_0' =
     \sqrt{r_0^2 + \rho_0^2} +s_0 - \sqrt{r_0^2 + ( \rho_0 + s_0)^2} \leq
      \frac{2r_0^2}{\sqrt{r_0^2+\rho_0^2} +\rho_0} \leq \epsilon$$

\vspace{-2mm}
\noindent for $t > T(r,\epsilon)$.
As $y'=\alpha (t')$ with $t'$ arbitrarily greater than $t$, taking the limit in (\ref{comparison}) for $t' \rightarrow \infty$   proves the lemma.\qed
 \vspace{4mm}

\vspace{2mm}
\emph{Proof of Proposition \ref{prophadamard}.}
\vspace{1mm}

\noindent {\em Let us first prove (iii).} Let $\alpha, \beta$ be rays with origins $a,b$ and initial conditions $u=\alpha'(0),v = \beta'(0)$ and let $K$ be any fixed compact set containing $a,b$. \linebreak 
We have to show that, for any arbitrary  $\delta>0$,  if  $u$ is sufficiently close to   $v$ then  $|B_{\alpha} (x,y) - B_{\beta} (x,y)| \!<\! \delta$ for all $x,y \in K$.
 Now, the Uniform Approximation Lemma ensures that we can replace   $B_{\alpha} (x,y)$ and  $ B_{\beta} (x,y)$ with $b_{\alpha (t)} (x,y)$ and $b_{\beta (t)}  (x,y)$,  making an error smaller than  $\delta/3$, by taking any $t> T(K,\delta/3)$. But the difference between these two functions is smaller than $2 d(\alpha (t),  \beta (t))$; and this, for any fixed $t$,  tends to zero as $u \rightarrow v$, on any Riemannian manifold. 
\vspace{2mm}

\noindent {\em  Let us now prove (ii).} \\
Assume that $(P_k)  \rightarrow \xi = B_{(P_k)} (o, \cdot)$. Then, consider the geodesic segments $\alpha_k= [o,P_k]$ and their velocity vector $u_k=\alpha_k'(0)$. Up to  a subsequence, the $u_k$'s converge to some unitary vector $u \in S_o \tilde{X}$.  As before, for any fixed compact set $K$,  the Uniform Approximation Lemma  ensures that 
$b_{\alpha_k(t)} (x, y) \simeq_{\epsilon} B_{\alpha_k} (x,y)$,  for any $t \geq T(K, \epsilon)$ and for all $x,y \in K$; in particular,
$b_{P_k} (x, y) \simeq_{\epsilon}  B_{\alpha_k} (x,y)$ if $t_k = d(o, P_k) > T(K, \epsilon)$. On the other hand, $B_{\alpha_k} (x,y) \simeq_{\epsilon} B_{\alpha} (x,y)$ if $k\gg 0$,  by (iii);  so passing to limits for $k \rightarrow \infty$, we deduce that $B_{(P_k)} (x,y)= B_{\alpha}(x,y)$ on $K$ and, as $K$ is arbitrary, $B_{(P_k)} = B_{\alpha}$.
\vspace{2mm}

\noindent  {\em We now prove   the first equivalence $d_{\infty} (\alpha, \beta) < \infty \;  \Leftrightarrow \; B_{\alpha}=B_{\beta}$ in (i)}. \\
Let $a=\alpha(0)$, $b=\beta(0)$ be the origins of $\alpha, \beta$.
If $d_{\infty} (\alpha, \beta) < \infty$,   by convexity of the distance in nonpositive curvature we deduce that there exist points $a_k, b_k$ tending to infinity  respectively along $\alpha$ and $\beta$, such that
\vspace{-3mm}

$$\lim_{k \rightarrow \infty}  a_k  b_k  = d = d(\alpha, \beta)$$

\vspace{-1mm}
\noindent Clearly,   the angles $\widehat{aa_kb_k}$ and  $\widehat{bb_ka_k}$ tend to $\pi/2$.
Now let $y$ be arbitrarily fixed, with $D= d(a,y)$.
By comparison with the Euclidean case, the tangent of the angle $\widehat{aa_ky}$ is  smaller than $D/(aa_k-D)$, which goes to zero for $k \rightarrow \infty$, so the angle  $\vartheta_k= \widehat{y a_k b_k} \rightarrow \pi/2$. Now we know, by comparison geometry,  that 
\vspace{-3mm}

$$ (b_k y) ^2 \geq  (a_k y) ^2 +  (a_k b_k) ^2 -2 a_k y \cdot  a_k b_k \cdot \cos  \vartheta_k$$

\vspace{-1mm}
\noindent hence $\liminf_{k \rightarrow \infty}  b_k y  - a_k y  \geq -\lim_k  a_k b_k \cos \vartheta_k =0$.
One proves analogously that $\liminf_{k \rightarrow \infty}  a_k y  -  b_k y  =0$, 
hence we deduce that  $\lim_{k \rightarrow \infty}   b_k  y  -  a_k  y  = 0$. As $y$ is arbitrary, this shows that $B_{\beta} = B_{\alpha}$.

\noindent Conversely, assume that $d_{\infty}(\alpha, \beta) = \infty$.
Up to possibly extending $\alpha$ and $\beta$ beyond their origins,   we may assume  that $a$ is the projection of $b$ over $\alpha$ and, {\em moreover,   that $\widehat{ab\beta(t)} \geq \pi/2$} (for $t \gg 0$).
In fact, let $\widetilde{\alpha}, \widetilde{\beta}$ be the bi-infinite \nolinebreak  geodesics extending $\alpha, \beta$:
either   $\limsup_{t \rightarrow -\infty}d(\widetilde{\alpha} (t),\widetilde{\beta}(t))$ is unbounded and, by convexity, there exists a minimal geodesic segment between $\widetilde{\alpha}$ and $\widetilde{\beta}$ (orthogonal to both  $\widetilde{\alpha},\widetilde{\beta}$);
or  $\limsup_{t \rightarrow -\infty}d(\widetilde{\alpha} (t), \widetilde{\beta}(t))$ is bounded, so the angle
$\widehat{\widetilde{\alpha}(t)a \widetilde{\beta}(t)} \;\; \rightarrow 0$  and   $[a,b,\widetilde{\beta}(t)]$ tends to the limit triangle  $\widetilde{\alpha} |_{\mathbb{R}^-}  \cup [a,b] \cup \widetilde{\beta} |_{\mathbb{R}^-} $ for $t\rightarrow -\infty$; as  the sum of its angles cannot exceed $\pi$, we deduce that   $\widehat{ab\beta(t) } \geq \frac{\pi}{2}$ for $t \gg0$.

\noindent So, now consider the triangle   $[a,b,\beta(t)]$ for $t \gg 0$. The angle $\widehat{\alpha(t)a\beta(t)}$ does not tend to zero for $t\rightarrow +\infty$, otherwise  $\alpha |_{\mathbb{R}^+} \cup  [a,b] \cup \beta |_{\mathbb{R}^+} $ would be again a limit triangle,  whose sum of angles necessarily would be $\pi$; thus, it would be  flat and totally geodesic, and $\lim_{t \rightarrow + \infty}d(\alpha(t), \beta(t)$ would be bounded.
Therefore,  $\widehat{\alpha(t)a\beta(t)}\geq \vartheta_0 >0$ for $t \rightarrow +\infty$.
By comparing  $[a,\alpha(s), \beta(t)]$, for $s,t \geq 0$, with an Euclidean triangle, we then get
\vspace{-3mm}

\begin{equation}
\label{carnot}
 (\alpha(s)  \beta (t))^2  \geq  s^2 +   (a \beta(t))^2 - 2 s \cdot a  \beta(t) \cdot  \cos \vartheta_0 
\end{equation}

\vspace{-1mm}
\noindent so $B_{\beta} (a, \alpha(s) ) = \lim_{t\rightarrow +\infty}   a \beta(t)  -  \beta (t) \alpha(s))  \leq s \cos \vartheta_0 < s = B_{\alpha} ( a, \alpha(s))$. \linebreak
This shows that $B_{\alpha} \neq B_{\beta}$.
\vspace{2mm}

\noindent  {\em  Proof of  the equivalence $B_{\alpha}=B_{\beta} \;  \Leftrightarrow \; \alpha \prec \beta$}. \\
One implication is true on any Riemannian manifold, as we have seen in Theorem {teorcoray}.
So, assume  that $\alpha \prec \beta$: let $\alpha_n=\overrightarrow{a_n b_n} \rightarrow \alpha$  with $a_n \rightarrow a$, $b_n=\beta(t_n)= \alpha_n(s_n)$ for $t_n, s_n \rightarrow +\infty$.
Let $K$ be a compact set containing $a,b$, the $a_n$ and points $x,y$,   and let $\epsilon>0$;
then, choose
$n \gg0 $ so that $s_n, t_n >T(K, \epsilon)$ of Lemma \ref{lemmaapprox} and such that $B_{\alpha_n} \eqsim_{\epsilon} B_{\alpha}$ on $K$, by (iii).
By Lemma \ref{lemmaapprox} and monotonicity of the Busemann cocycle we then get 
\vspace{-6mm}

$$B_{\alpha}(x,y) \eqsim_{\epsilon} B_{\alpha_n}(x,y)
\eqsim_{\epsilon} b_{\alpha_n (s_n)} (x, y)
= b_{\beta_n (t_n)} (x, y)
\eqsim_{\epsilon} B_{\beta} (x,y)$$

\vspace{-2mm}
\noindent and as $\epsilon$ is arbitrary, we deduce that $B_{\alpha}(x,y)=B_{\beta} (x,y)$.
\vspace{2mm}

\noindent {\em Finally},  if two rays $\alpha, \beta$ with common origin $o$  make angle  $\vartheta_0 \neq 0$, then $d(\alpha(s), \beta(t))$ grows at least as in the Euclidean case according to the formula (\ref{carnot}), hence the rays are not Busemann equivalent, so the restriction of the Busemann map ${\cal R}_o (X) \rightarrow \partial X$ is injective.\qed

\subsection{Hyperbolic computations.}
\label{apphyp}
 
\noindent {\bf Lemma \ref{lemmahyp1}}
{\em Let $\tilde{o} \in \Hy^2$, and let $C, C'$ two ultraparallel geodesics (i.e. with no common point in $\Hy^2 \cup \partial \Hy^2$) such that $d(\tilde{o} , C)=d(\tilde{o} , C')$. Then:

\noindent (i) there exists a unique hyperbolic isometry $g$ with axis perpendicular to $C,C'$ and
 such that $g(C)=C'$;

\noindent (ii) $g^{-1}\tilde{o}$ and $g\tilde{o}$ are obtained, respectively, by the hyperbolic reflections of $\tilde{o}$ with respect to $C,C'$;

\noindent (iii) the Dirichlet domain $D(g,\tilde{o})$ has boundary $C \cup C'$.
}
 \vspace{2mm}

 {\em Proof.} By convexity of the distance function, there exists a unique common perpendicular $\tilde{g}$ to  $C, C'$, so $g$ is the unique hyperbolic translation along $\tilde{g}$ sending $C$ to $C'$. Let $\Delta (g)$ the displacement of $g$, let $\tilde{o}_0$ be the projection of $\tilde{o}$ on $\tilde{g}$, and  let $p=C\cap \tilde{g}$.  By symmetry,  $\Delta(g) = d(C,C')= 2 \tilde{o}_0 p$. Now consider the hyperbolic reflection $R$ with respect to $C$, and define 
  $\tilde{c}=R(\tilde{o})$,  $\tilde{c}_0= R(\tilde{o}_0)$ and $q=[\tilde{o},\tilde{c}] \cap C$. Since $\tilde{g}$ is perpendicular to $C$, $R$ preserves $\tilde{g}$; we deduce that  $ [\tilde{c},\tilde{c}_0]=R([\tilde{o},\tilde{o}_0])$ is also perpendicular to $\tilde{g}$. As $\tilde{o}_0\tilde{c}_0=2\tilde{o}_0p=\Delta(g) $, it follows that $g^{-1}\tilde{o}=\tilde{c}$. Then, $C$ is one of the two boundaries of $D(g,\tilde{o})$, as it is the perpendicular bisector of $[\tilde{o},\tilde{c}]$. The verification for $g\tilde{o}$ and  $C'$  is the same.\qed

\begin{lemma}
\label{lemmalength}
\hspace{-2mm} There exists a positive function $\delta(t, \epsilon)$ for $t,\! \epsilon \!>\! 0$, increasing in \nolinebreak $t$, with the following property. Let $\alpha$ be any geodesic of $\Hy^2$  and assume that  $p_1, p_2$ are points  with  projections  $q_1,q_2$ on $\alpha$ such that  $d(q_1,q_2)=t$: if  $d(p_1, \alpha) =\epsilon$, then $d(p_1, p_2) \geq t + \delta(t, \epsilon)$.
\end{lemma}

\vspace{2mm}
{\em Proof.} Consider the projection $p_1'$ of $p_1$ on the geodesic containing $p_2,q_2$ and let $d=p_1 p'_1 \leq p_1 p_2$.  Let $c=p_1q_2$ and $\beta = \widehat{p_1 q_2 q_1}$.
By the sinus and cosinus formula applied, respectively, to the triangles $[p_1,p'_1,q_2]$ and $[p_1,q_1,q_2]$ we find
\vspace{-3mm}

$$\sinh d = \sinh c \cdot \cos \beta = \cosh c  \cdot \tanh t$$

\vspace{-1mm}
\noindent and by Phytagora's formula we deduce that $\sinh d =  \cosh \epsilon \sinh t$. This shows that $d = t + \delta (t, \epsilon)$, for a  positive function $\delta (t, \epsilon)$ when $t, \epsilon >0$.
To see that $\delta (t, \epsilon)$ is increasing with $t$ we just compute the derivative 
\vspace{-3mm}

$$\partial_t \, \delta(t,\epsilon) = d(t)' -1 = \frac{\cosh \epsilon \cosh t }{\cosh d} - 1 = \frac{\cosh c}{\cosh d} - 1 >0$$

\vspace{-1mm}
\noindent as $c>d$ for $\epsilon >0$.\qed

\pagebreak

\end{document}